\documentclass[11pt]{amsart}
\usepackage{amssymb,amsmath,epsfig,graphics}
\theoremstyle{plain}
\newtheorem{thrm}{Theorem}[section]
\newtheorem{lemma}[thrm]{Lemma}
\newtheorem{prop}[thrm]{Proposition}
\newtheorem{cor}[thrm]{Corollary}
\newtheorem{rmrk}[thrm]{Remark}
\newtheorem{dfn}[thrm]{Definition}

\numberwithin{equation}{section}
\numberwithin{figure}{section}

\begin{document}
\newcommand{\SL}{\mathcal L^{1,p}(\Om)}
\newcommand{\Lp}{L^p(\Omega)}
\newcommand{\CO}{C^\infty_0(\Omega)}
\newcommand{\Rn}{\mathbb R^n}
\newcommand{\Rm}{\mathbb R^m}
\newcommand{\R}{\mathbb R}
\newcommand{\Om}{\Omega}
\newcommand{\Hn}{\mathbb H^n}
\newcommand{\HH}{\mathbb H^1}
\newcommand{\eps}{\epsilon}
\newcommand{\BVX}{BV_H(\Omega)}
\newcommand{\IO}{\int_\Omega}
\newcommand{\bG}{\boldsymbol{G}}
\newcommand{\bg}{\mathfrak g}
\newcommand{\p}{\partial}
\newcommand{\Xnu}{\overset{\rightarrow}{ H_\nu}}
\newcommand{\nuX}{\boldsymbol{\nu}_H}
\newcommand{\Up}{\boldsymbol{\mathcal Y}_H}
\newcommand{\n}{\boldsymbol \nu}
\newcommand{\sigmau}{\boldsymbol{\sigma}^u_H}
\newcommand{\di}{\delta_{H,i}}
\newcommand{\del}{\delta_H}
\newcommand{\nui}{\nu_{H,i}}
\newcommand{\nuj}{\nu_{H,j}}
\newcommand{\dej}{\delta_{H,j}}
\newcommand{\cx}{\boldsymbol{c}_S}
\newcommand{\sx}{\sigma_H}
\newcommand{\lx}{\mathcal L_H}
\newcommand{\rad}{\overline u}
\newcommand{\nH}{\nabla_H}
\newcommand{\pb}{\overline p}
\newcommand{\qb}{\overline q}
\newcommand{\nuu}{\boldsymbol \nu_{H,u}}
\newcommand{\nuv}{\boldsymbol \nu_{H,v}}
\newcommand{\Bl}{\Bigl|_{\lambda = 0}}
\newcommand{\AD}{\mathcal A_{\mathcal D}}
\newcommand{\mS}{\mathcal S}
\newcommand{\delh}{\Delta_H}
\newcommand{\delinf}{\Delta_{H,\infty}}
\newcommand{\nabh}{\nabla_H}
\newcommand{\delp}{\Delta_{H,p}}
\newcommand{\mO}{\mathcal O}
\newcommand{\delhs}{\Delta_{H,S}}
\newcommand{\lhs}{\hat{\Delta}_{H,S}}
\newcommand{\bN}{\boldsymbol{N}}
\newcommand{\bnu}{\boldsymbol \nu}
\newcommand{\la}{\lambda}
\newcommand{\nup}{\boldsymbol{\nu}_H^\perp}


\title[A partial solution of the isoperimetric problem for the Heisenberg group]
{A partial solution of the isoperimetric problem for the Heisenberg group}

\author{Donatella Danielli}
\address{Department of Mathematics \\
Purdue University\\
West Lafayette, IN 47907}
\email[Donatella Danielli]{danielli@math.purdue.edu}
\thanks{First author supported in part by NSF grants DMS-0002801 and CAREER DMS-0239771}
\author{Nicola Garofalo}
\address{Department of Mathematics\\Purdue University \\
West Lafayette, IN 47907} \email[Nicola
Garofalo]{garofalo@math.purdue.edu}
\address{Dip. Metodi e Modelli Matematici per le Scienze
Applicate\\
 Univ. Padova, 35100 Padova, Italy}
\email[Nicola Garofalo]{garofalo@math.purdue.edu}
\thanks{Second author supported in part by NSF Grants DMS-0070492 and DMS-0300477}
\author{Duy-Minh Nhieu}
\address{Department of Mathematics \\
Georgetown University \\
Washington DC 20057-1233} \email[Duy-Minh Nhieu]{dn9@georgetown.edu}
%
%

\maketitle

\baselineskip 16pt

\tableofcontents


\section{\textbf{Introduction}}

\vskip 0.2in

The classical isoperimetric problem states that among all
measurable sets with assigned volume the ball minimizes the
perimeter. This is the content of the celebrated isoperimetric
inequality, see \cite{DG3},
\begin{equation}\label{ii}
|E|^{\frac{n-1}{n}}\ \leq\ C_n\ P(E)\ ,
\end{equation}
which holds for all measurable sets $E\subset \Rn$ with constant
$C_n = n \sqrt \pi/\Gamma(n/2 + 1)^{1/n}$. In \eqref{ii}, $P(E)$
denotes the perimeter in the sense of De Giorgi, see \cite{DG1},
\cite{DG2}, i.e., the total variation of the indicator function of
$E$. Equality holds in \eqref{ii} if and only if (up to negligible
sets) $E = B(x,R) = \{y\in \Rn\mid |y-x|<~R\}$,
 a Euclidean ball. It is well-known that \eqref{ii} is
 equivalent to the geometric Sobolev inequality for $BV$ functions,
 see \cite{FR}. An analogous ``isoperimetric inequality" was proved
 in \cite{GN} in the general setting of a Carnot-Carath\'eodory
 space, and such inequality was used, among other things, to
 establish a geometric embedding for horizontal $BV$ functions,
 similar to Fleming and Rishel's one. However, the question of the
 optimal configurations in such isoperimetric inequality was left
 open.

The aim of this paper is to bring a partial solution to this open
problem in the Heisenberg group $\Hn$. We recall that $\Hn$ is the
simplest and perhaps most important prototype of a class of
nilpotent Lie groups, called Carnot groups, which play a
fundamental role in analysis and geometry, see \cite{Ca},
\cite{Ch}, \cite{H}, \cite{St}, \cite{Be}, \cite{Gro1},
\cite{Gro2}, \cite{E1}, \cite{E2}, \cite{E3}, \cite{DGN2}. Its
underlying manifold is $\mathbb R^{2n+1}$ with non-commutative
group law
\begin{equation}\label{Hn}
g\ g'\ =\ (x,y,t)\ (x',y',t')\ =\ (x + x', y + y', t + t' +
\frac{1}{2} ( <x,y'> - <x',y> ))\ ,
\end{equation}
where we have let $x,x', y, y' \in \Rn$, $t,t'\in \mathbb R$. If
$L_g(g') = g g'$ denotes the operator of left-translation, let
$(L_g)_*$ indicate its differential. The Heisenberg algebra admits
the decomposition $\mathfrak h_n = V_1\oplus V_2$, where $V_1 =
\mathbb R^{2n} \times \{0\}$, and $V_2 = \{0\} \times \R$.
Identifying $\mathfrak h_n$ with the space of left-invariant
vector fields on $\Hn$, one easily recognizes that a basis for
$\mathfrak h_n$ is given by the $2n+1$ vector fields
\begin{equation}\label{vf2}
\begin{cases}
(L_g)_*\left(\frac{\p}{\p x_i}\right)\ \overset{def}{=}\ X_i\ =\
\frac{\p}{\p x_i} - \frac{y_i}{2}\ \frac{\p}{\p t}\ ,
\\
(L_g)_*\left(\frac{\p}{\p y_i}\right)\ \overset{def}{=}\ X_{n+i}\
=\ \frac{\p}{\p y_i} + \frac{x_i}{2}\ \frac{\p}{\p t}\ ,
\\
(L_g)_*\left(\frac{\p}{\p t}\right)\ \overset{def}{=}\ T\ =\
\frac{\p}{\p t}\ ,
\end{cases}
\end{equation}
and that the only non-trivial commutation relation is
\begin{equation}\label{commHn}
[X_i,X_{n+j}]\ =\ T\ \delta_{ij}\ ,\quad\quad\quad\quad i , j =
1,...,n\ .
\end{equation}

In \eqref{vf2} we have identified the standard basis
$\{e_1,...,e_{2n}, e_{2n + 1}\}$ of $\mathbb R^{2n+1}$ with the
system of (constant) vector fields $\{\p/\p x_1,...,\p/\p
y_{n},\p/\p t\}$. Because of \eqref{commHn} we have $[V_1,V_1] =
V_2$, $[V_1,V_2] = \{0\}$, thus $\Hn$ is a graded nilpotent Lie
group of step $r=2$. Lebesgue measure $dg = dzdt$ is a bi-invariant
Haar measure on $\Hn$. If we denote by $\delta_\lambda(z,t) =
(\lambda z, \lambda^2 t)$ the non-isotropic dilations associated
with the grading of the Lie algebra, then $d(\delta_\lambda g) =
\lambda^Q dg$, where $Q = 2n + 2$ is the homogeneous dimension of
$\Hn$.

In what follows we denote by $P_H(E;\Hn)$ the intrinsic, or
$H$-perimeter of $E\subset \Hn$ associated with the
bracket-generating system $X=\{X_1,...,X_{2n}\}$. Such notion will
be recalled in Section \ref{S:isoperimetric}. To state our theorem
we let $\Hn_+ = \{(z,t) \in \Hn \mid t
> 0\}$, $\Hn_- = \{(z,t) \in \Hn \mid t < 0\}$, and consider the
collection
\[
\mathcal E\ =\ \{E\subset \Hn\,|\, E \quad\text{satisfies
$(i)-(iii)$}\}\ , \]
 where
\begin{itemize}
\item[(i)] $|E\cap \mathbb H^n_{+}|\ =\ |E \cap \mathbb H^n_{-}|$ ;
\item[(ii)] there exist $R>0$, and functions $u,v:  \overline B(0,R)\to[0,\infty)$, with $u , v \in C^2(B(0,R)) \cap C(\overline B(0,R))$, $u = v =
0$ on $\p B(0,R)$, and such that
\[
\partial E \cap \mathbb
H^n_{+}\ =\ \{(z,t)\in \Hn_+\,|\,\ |z|<R \ , \ t\ =\ u(z)\}\ ,
\]
\[
 \partial E \cap \mathbb H^n_{-}\ =\ \{(z,t)\in \Hn_{-}\,|\,\ |z|<R\ ,\ t\ =\ -\ v(z)\}\ .
 \]
\item[(iii)] $\{z \in B(0,R)\mid u(z) = 0\}\ \cap \ \{z\in B(0,R)\mid v(z) =
 0\}\ =\ \varnothing$\ .
\end{itemize}

\begin{figure}[h]
  \hskip 0.3in
  \includegraphics{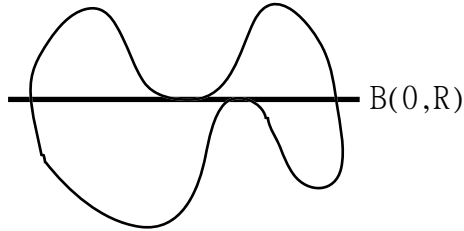}
  \caption{$E\in \mathcal E$}
  \label{pic:E}
\end{figure}

We note explicitly that condition (iii) serves to guarantee that
every $E\in \mathcal E$ is a piecewise $C^2$ domain in $\Hn$ (with
possible discontinuities in the derivatives only on that part of $E$
which intersects the hyperplane $t=0$). We also stress that the
upper and lower portions of a set $E\in \mathcal E$ can be described
by possibly different $C^2$ graphs, and that, besides $C^2$
smoothness, and the fact that their common domain is a ball, no
additional assumption is made on the functions $u$ and $v$. For
instance, we do not require a priori that $u$ and/or  $v$ are
spherically symmetric. Here is our main result.

\medskip

\begin{thrm}\label{T:isoprofile}
Let $V>0$, and define the number $R>0$ by
\[ R\ =\ \left(\frac{(Q-2) \Gamma\left(\frac{Q+2}{2}\right)
\Gamma\left(\frac{Q-2}{2}\right)}{\pi^{\frac{Q-1}{2}}
\Gamma\left(\frac{Q+1}{2}\right)}\right)^{1/Q}\ V^{1/Q}\ .
\]
Given such $R$, then the variational problem
\[
\underset{E \in \mathcal E , |E|=V}{\min}\ P_H(E;\Hn)
\]
has a unique solution $E_R = \delta_R(E_o) \in\mathcal E$, where
$\partial E_o$ is described by the graph $t = \pm\ u_o(z)$, with
\begin{align}\label{iso}
 u_o(z)\ \overset{def}=\   \left\{\frac{\pi}{8}\ +\
\frac{|z|}{4} \sqrt{1-|z|^2}\ -\ \frac{1}{4}\,\sin^{-1}(|z|)
\right\}\ , \quad |z|\ \leq\ 1\ .
\end{align}
The sign $\pm$ depends on whether one considers $\partial E_o \cap
\Hn_+$, or $\partial E_o\cap \Hn_-$. Finally, the boundary
$\partial E_R = \delta_R(\p E_o)$ of the bounded open set $E_R$ is
only of class $C^2$, but not of class $C^3$, near its two
characteristic points $\left(0,\pm \frac{\pi R^2}{8}\right)$, it
is $C^\infty$ away from them, and $S_R = \p E_R$ has positive
constant $H$-mean curvature given by
\[
\mathcal H\ =\ \frac{Q - 2}{R}\ .
\]
\end{thrm}

\medskip

\begin{rmrk}\label{R:mm}
We notice explicitly that the function $u_o$ in \eqref{iso} can
also be expressed as follows \[ u_o(z)\ =\ \frac{1}{2}
\int_{\sin^{-1}(|z|)}^{\frac{\pi}{2}} \sin^2 \tau\ d\tau\ . \]
\end{rmrk}

\medskip

\begin{rmrk}\label{R:smooth}
We emphasize that, as the reader will recognize, for our proof of
the existence of a global minimizer it suffices to assume that the
two functions $u$ and $v$ in the definition of the sets of the class
$\mathcal E$ are $C^{1,1}_{loc}(B(0,R))$. It is an open question
whether $u,v\in C^1(B(0,R))$ is enough. This is possible thanks to a
sharp result of Balogh concerning the size of the characteristic
set, see Theorem \ref{T:balogh} below. In our proof of the
uniqueness of the global minimizer, instead, it is convenient to
work under the hypothesis of $C^2$ smoothness. However, with little
extra care, it should be possible to relax it to $C^{1,1}_{loc}$.
\end{rmrk}

\medskip

For the notion of $H$-mean curvature of a $C^2$ hypersurface $\mS
\subset \Hn$ we refer the reader to Definition \ref{D:HMC} in
Section \ref{S:PS}. This notion of horizontal mean curvature, which
is of course central to the present study, was introduced in
\cite{DGN3}. Its geometric interpretation is that, in the
neighborhood of a non-characteristic point $g\in \mS$, it coincides
with the standard Riemannian mean curvature of the
$2n-1$-dimensional submersed manifold obtained by intersecting the
hypersurface $\mS$ with the fiber of the horizontal subbundle
$H_g\Hn$, see also \cite{DGN4} where a related notion of Gaussian
curvature was introduced. A seemingly different notion, based on the
Riemannian regularization of the sub-Riemannian metric of $\Hn$, was
proposed in \cite{Pa}, but the two are in fact equivalent, see
\cite{DGN3}. From Theorem \ref{T:isoprofile} we obtain the following
isoperimetric inequality for the horizontal perimeter.

\medskip

\begin{thrm}\label{T:isoine}
Let $\mathcal E$ be as above, and denote by $\tilde{\mathcal E}$
the class of sets of the type $\delta_\lambda L_g(E)$, for some
$E\in \mathcal E$, $\lambda >0$ and $g\in \Hn$, then the following
isoperimetric inequality holds
\begin{equation}\label{isoine}
|E|^{\frac{Q-1}{Q}}\ \leq\ C_Q\ P_H(E;\Hn)\ ,\quad\quad\quad\quad E
\in \tilde{ \mathcal E}\ ,
\end{equation}
where
\[
C_Q\ =\ \frac{(Q-1)\Gamma\left(\frac{Q}{2}\right)^\frac{2}{Q}}
{Q^\frac{Q-1}{Q}(Q-2)
\Gamma\left(\frac{Q+1}{2}\right)^\frac{1}{Q}\pi^\frac{Q-1}{2Q}}\ ,
\]
with equality if and only if for some $\lambda >0$ and $g\in \Hn$
one has $E = L_g \delta_\lambda (E_o)$, where $E_o$ is given by
\eqref{iso}.
\end{thrm}

\medskip

Fig.1.1 gives a representation of the isoperimetric set $E_o$ in
Theorem \ref{T:isoprofile} in the special case $n=1$. We note that
the invariance of the isoperimetric quotient with respect to the
group left-translations $L_g$ and dilations $\delta_\la$ is
guaranteed by Propositions \ref{P:invariance} and
\ref{P:invariance2}.

\begin{figure}[h]
  \hskip 0.3in
  \includegraphics{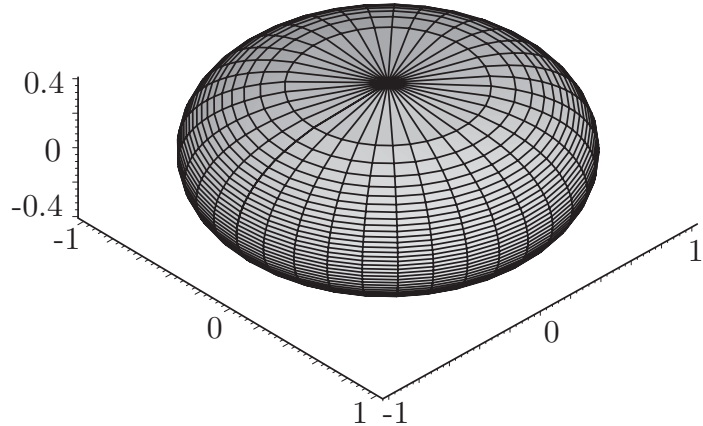}
  \caption{Isoperimetric set in $\mathbb H^1$ with $R=1$}
  \label{pic:Iso_set}
\end{figure}

\noindent A remarkable property of the isoperimetric sets is that,
similarly to their Riemannian predecessors, they have constant
$H$-mean curvature. It is tempting, and also natural, to
conjecture that the set $E_o$ described by \eqref{iso}, along with
its left-translated and dilated, exaust \emph{all} the
isoperimetric sets in $\Hn$ (for the definition of such sets, see
Definition \ref{D:isoquotient} below). By this we mean that
Theorem \ref{T:isoine} continues to be valid when one replaces the
class $\tilde{\mathcal E}$ with that of all measurable sets
$E\subset \Hn$ with locally finite $H$-perimeter. At the moment,
this remains a challenging open problem. In this connection,
another interesting conjecture is as follows: \emph{Let $\mathcal
S\subset \Hn$ be a $C^2$, compact oriented hypersurface. Suppose
that for some $\alpha >0$
\begin{equation}\label{conj2}
\mathcal H\ \equiv\ \alpha \ \quad\quad\quad\text{on}\quad
\mathcal S\ .
\end{equation}
Is it true that, up to a left translation, if we denote by $\mathcal
S^+ = \mS \cap \Hn_+$, $\mathcal S^- = \mS \cap \Hn_-$, then
$\mathcal S^+$, $\mathcal S^-$ are respectively described by
\begin{align}\label{isoR}
t\ =\ & \pm\ \left\{\frac{1}{4}|z|\sqrt{R^2-|z|^2}\ -\
\frac{R^2}{4}\,\tan^{-1}\left(\frac{|z|}{\sqrt{R^2-|z|^2}}\right)\
+\ \frac{\pi R^2}{8}\right\}\ , \quad |z|\ \leq\ R\ ,
\end{align}
where $R = (Q-2)/\alpha$?} Concerning this conjecture we remark that
Theorem \ref{T:isoprofile} provides evidence in favor of it. As it
is well-known, the Euclidean counterpart of it is contained in the
celebrated \emph{soap bubble} theorem of A.D. Alexandrov \cite{A}.
We mention that, after this paper was completed, we have received an
interesting preprint from Ritor\'e and Rosales \cite{RR2} in which,
among other results, the authors prove the above soap bubble
conjecture in the first Heisenberg group $\HH$.

To put the above results in a broader perspective we recall that
in any Carnot group a general scale invariant isoperimetric
inequality is available. In fact, using the results in \cite{CDG},
\cite{GN} one can prove the following  theorem, see Theorem
\ref{T:relisoG} in Section \ref{S:isoperimetric}.

\medskip

\begin{thrm}\label{T:relisoG0}
Let $\bG$ be a Carnot group with homogeneous dimension $Q$. There
exists a constant $C_{iso}(\bG)>0$ such that, for every
$H$-Caccioppoli set $E\subset \bG$, one has
\[
|E|^{(Q-1)/Q}\ \leq\ C_{iso}(\bG)\ P_H(E;\bG)\ .
\]
\end{thrm}

\medskip

A measurable set $E\subset \bG$ is called a $H$-Caccioppoli set if
$P_H(E;\omega)<\infty$ for any $\omega \subset \subset \bG$.
Theorem \ref{T:relisoG0} generalizes an earlier result of Pansu
\cite{P}, who proved a related inequality for the first Heisenberg
group $\mathbb H^1$, but with the $H$-perimeter in the right-hand
side replaced by the $3$-dimensional Hausdorff measure $\mathcal
H^3$ in $\mathbb H^1$ constructed with the Carnot-Carath\'eodory
distance associated with the horizontal subbundle $H\mathbb H^1$
defined by $\{X_1,X_2\}$ in \eqref{vf2}. One should keep in mind
that the homogeneous dimension of $\mathbb H^1$ is $Q = 4$, so $3=
Q-1$, which explains the appearance of $\mathcal H^3$ in Pansu's
result. It should also be said that some authors attribute to
Pansu \cite{P2} the conjecture that the isoperimetric sets in
$\mathbb H^1$ have the form \eqref{iso}. We mention that other
isoperimetric and Fleming-Rishel type Gagliardo-Nirenberg
inequalities have been obtained by several authors at several
times, see \cite{Va1}, \cite{Va2}, \cite{VSC}, \cite{CS},
\cite{BM}, \cite{FGW}, \cite{MaSC}. We now introduce the following
definition.

\medskip

\begin{dfn}\label{D:isoquotient}
Given a Carnot group $\bG$ with homogeneous dimension $Q$ we
define the \emph{isoperimetric constant} of $\bG$ as
\[
\alpha_{iso}(\bG)\ =\ \underset{E\subset \bG}{\inf}\
\frac{P_H(E;\bG)}{|E|^{(Q-1)/Q}}\ ,
\]
where the infimum is taken on all $H$-Caccioppoli sets $E$ such
that $0 < |E| < \infty$. If a measurable set $E_o$ is such that
\[
\alpha_{iso}(\bG)\ =\ \frac{P_H(E_o;\bG)}{|E_o|^{(Q-1)/Q}}\ ,
\]
then we call it an isoperimetric set in $\bG$.
\end{dfn}

\medskip

We stress that, thanks to Theorem \ref{T:relisoG0}, the
isoperimetric constant is strictly positive. It should also be
observed that, using the representation formula for the
$H$-perimeter \begin{equation}\label{permeasure}
 P_H(E;\bG)\ =\
\int_{\p E} \frac{W}{|\bN|}\ dH_{N-1}\ , \end{equation} valid for
any bounded open set $E\subset \bG$ of class $C^1$, with Riemannian
outer normal $\bN$ and angle function $W = \sqrt{p_1^2 + ... +
p_m^2}$ (see Lemma \ref{L:smoothsets}, and \eqref{pi}, \eqref{W}),
one immediately recognizes that, since for any $\omega \subset
\subset \bG$ one has $W \leq C(\omega) |\bN|$, then $P_H(E;\bG) \leq
C\ H_{N-1}(\partial E) < \infty$. As a consequence,
$\alpha_{iso}(\bG) < \infty$ as well. What is not obvious instead is
the existence of isoperimetric sets. In this regard, one has the
following basic result proved in \cite{LR}.

\medskip

\begin{thrm}\label{T:existence0}
Let $\bG$ be a Carnot group, then there exists a bounded
$H$-Caccioppoli set $F_o$ such that
\[
P_H(F_o;\bG)\ =\ \alpha_{iso}(\bG)\ |F_o|^{(Q-1)/Q}\ .
\]
The equality continues to be valid if one replaces $F_o$ by
$L_{g_o} \circ \delta_\lambda (F_o)$, for any $\lambda > 0$,
$g_o\in \bG$.
\end{thrm}

\medskip

Of course, this result leaves open the fundamental question of the
classification of such sets. We stress that, in the generality of
Theorem \ref{T:relisoG0}, this problem is presently totally out of
reach. When $\bG = \Hn$, however, Theorems \ref{T:isoprofile} and
\ref{T:isoine} provide some basic progress in this direction. Our
main contribution is to use direct methods of the Calculus of
Variations to prove that the critical point \eqref{isoR} is a global
minimizer in the class $\mathcal E$. Furthermore, such global
minimizer is unique (modulo left-translations and dilations) in such
class. These results follow from some delicate properties of
convexity, and strict convexity at the global minimizer, of the
$H$-perimeter functional subject to a volume constraint.

In connection with our work, we mention that several authors have
recently studied the isoperimetric problem in $\Hn$, but under the
restriction that the class of competitors be $C^2$ smooth and
cylindrically symmetric, i.e., spherical symmetry about the $t$-axis
of the graph of the competing sets. For instance, in the recent
interesting work \cite{BC}, for the first Heisenberg group $\mathbb
H^1$, the authors prove that the flow by $H$-mean curvature of a
$C^2$ surface which is convex, and which is described by $t = \pm
f(|z|)$, with $f'<0$, converges to the sets \eqref{iso}. Notice,
however, that $f$ is spherically symmetric, convex, and that it is
assumed that the upper and lower part of the surface are described
by the same strictly decreasing function $f$. We also mention the
paper \cite{Pa} in which the author, still for $\mathbb H^1$,
heuristically derives the surface described by \eqref{iso} by
imposing the condition of constant $H$-mean curvature among all
$C^2$ surfaces which can be described by $t = \pm f(|z|)$. Recently,
Hladky and Pauls in \cite{HP} have proposed a general geometric
framework, which they call Vertically Rigid manifolds, and which
encompasses the class of Carnot groups, in which they study the
isoperimetric and the minimal surface problems. In this setting they
introduce a notion of horizontal mean curvature, and they show, in
particular, that remarkably the isoperimetric sets have constant
horizontal mean curvature. In the paper \cite{LM} the authors prove,
among other interesting results, that the $u_o$ in our Theorem
\ref{T:isoprofile} is a critical point (but not the unique global
minimizer) of the $H$-perimeter, when the class of competitors is
restricted to $C^2$ domains, with defining function of the type $t =
\pm f(|z|)$. A similar result has been also obtained in the
interesting recent preprint \cite{RR}, which also contains a
classification of the Delaunay type surfaces in $\Hn$. In this
connection, we also mention the earlier paper \cite{To}, in which
the author describes the Delaunay type surfaces of revolution in
$\HH$, heuristically computes the special solutions \eqref{iso}, and
shows that standard Schwarz symmetrization does not work in the
Heisenberg group. In \cite{FMP} the authors gave a complete
classification of the constant mean curvature surfaces (including
minimal) which are invariant with respect to $1$-dimensional closed
subgroups of ${\rm Iso}_0(H_3,g)$. We also mention the paper
\cite{Mo1}, in which the author proved that the
Carnot-Carath\'eodory ball in $\Hn$ is not an isoperimetric set.
Subsequently, in \cite{Mo2} he proved that, as a consequence of this
fact, a generalization of the Brunn-Minkowski inequality to $\Hn$
fails. Finally, in their interesting paper \cite{MoM} the authors
have established the isoperimetric inequality for the
Baouendi-Grushin vector fields $X_1 = \p_x$, $X_2 = |x|^\alpha
\p_t$, $\alpha >0$, in the plane $(x,t)$, and explicitly computed
the isoperimetric profiles. In the special case $\alpha = 1$, such
profiles are identical (up to a normalization of the vector fields)
to our $u_o$ in Theorem \ref{T:isoprofile}, see Remark \ref{R:mm}
above.

\medskip

\noindent \textbf{Acknowledgment$^1$:} For the first Heisenberg
group $\mathbb H^1$, and under the assumption that the isoperimetric
profile be of class $C^2$ and of the type $t = f(|z|)$, the idea of
using calculus of variations to explicitly determine $f(|z|)$, first
came about in computations that Giorgio Talenti and the second named
author carried in a set of unpublished notes in Oberwolfach in 1995.
We would like to thank G. Talenti for his initial contribution to
the present study. \footnotetext[1]{The results in this paper were
presented by the second named author in the lecture: "Remarks on the
best constant in the isoperimetric inequality for the Heisenberg
group and surfaces of constant mean curvature", Analysis seminar,
University of Arkansas, April 12, 2001,
(http://comp.uark.edu/$\sim$lanzani/schedule.html), by the third
named author at the international meeting on ``Subelliptic equations
and sub-Riemannian geometry", Arkansas, March 2003, and by the first
named author in the lecture ``Hypersurfaces of minimal type in
sub-Riemannian geometry", Seventh New Mexico Analysis Seminar,
University of New Mexico, October 2004.}

\vskip 0.6in


\section{\textbf{Isoperimetric inequalities in Carnot groups}}\label{S:isoperimetric}

\vskip 0.2in

The appropriateness of the notion of $H$-perimeter in
Carnot-Carath\'eodory geometry is witnessed by the isoperimetric
inequalities. Similarly to their Euclidean counterpart, these
inequalities play a fundamental role in the development of
geometric measure theory. Theorem \ref{T:relisoG0} represents a
sub-Riemannian analogue of the classical global isoperimetric
inequality. Such result can be extracted from the isoperimetric
inequalities obtained  in \cite{CDG} and \cite{GN}, but it is not
explicitly stated in either paper. Since a proof of Theorem
\ref{T:relisoG0} is not readily available in the literature, for
completeness we present it in this section.

Given a Carnot group $\bG$, its Lie algebra $\bg$ satisfies the
properties $\bg = V_1\oplus ... \oplus V_r$, where $[V_1,V_j] =
V_{j+1}$, $j=1,...,r-1$, and $[V_1,V_r] = \{0\}$. If $m_j = dim\
V_j$, $j=1,...,r$, then the homogeneous dimension of $\bG$ is
defined by $Q = m_1 + 2 m_2 + ... + r m_r$. The non-isotropic
dilations associated with the grading of $\bg$ are given by
$\Delta_\lambda(\xi_1 + ... + \xi_r) = \lambda \xi_1 + ... +
\lambda^r \xi_r$. Via the exponential mapping $\exp : \bg \to \bG$,
which is a global diffeomorphism onto, such dilations induce a
one-parameter group of dilations on $\bG$ as follows
$\delta_\lambda(g) = \exp \circ \Delta_\lambda \circ \exp^{-1} (g)$.
The push forward through $\exp$ of the standard Lebesgue measure on
$\bg$ is a bi-invariant Haar measure on $\bG$. We will denote it by
$dg$. Clearly, $d(\delta_\lambda g) = \lambda^Q\ dg$. For
simplicity, we let $m = m_1$. We fix orthonormal basis $\{e_1, ... ,
e_m\}$, .... , $\{e_{r,1}, ... , e_{r,m_r}$\}, of the layers $V_1,
... , V_r$, and consider the corresponding left-invariant vector
fields on $\bG$ defined by $X_1(g) = (L_g)_*(e_1), ... , X_m(g) =
(L_g)_*(e_m)$, ... , $X_{r,1}(g) = (L_g)_*(e_{r,1}), ... ,
X_{r,m_r}(g) = (L_g)_*(e_{r,m_r})$. We will assume that $\bG$ is
endowed with a left-invariant Riemannian metric $<\cdot,\cdot>$ with
respect to which these vector fields constitute and orthonormal
basis. No other inner product will be used in this paper. We denote
by $H\bG\subset T\bG$ the subbundle of the tangent bundle generated
by $\{X_1,...,X_m\}$. We next recall the notion of $H$-perimeter,
see e.g. \cite{CDG}. Given an open set $\Om\subset \bG$, we let
\[
\mathcal F(\Om)\ =\ \left\{\zeta = \sum_{i=1}^m \zeta_i X_i \in
\Gamma^1_0(\Om,H\bG)\ \mid\ |\zeta|_\infty\ =\ \sup_{\Om}\
|\zeta|\ =\ \sup_{\Om} \left(\sum_{i=1}^m \zeta_i^2\right)^{1/2}
\leq 1\right\}\ ,
\]
where we say that $\zeta \in \Gamma^1_0(\Om,H\bG)$ if $X_j \zeta_i
\in C_0(\Om)$ for $i,j = 1,...,m$. Given $\zeta \in
\Gamma^1_0(\Om,H\bG)$ we define
\[ div_H\ \zeta\ =\
 \sum_{i=1}^m X_i \zeta_i\ .
 \]

 For a function $u\in L^1_{loc}(\Om)$, the $H$-variation of $u$
with respect to $\Om$ is defined by
\[
Var_H(u;\Om)\ =\ \underset{\zeta\in \mathcal F(\Om)}{\sup}\
\int_{\bG} u\ div_H \zeta\ dg\ .
\]

We say that $u\in L^1(\Om)$ has bounded $H$-variation in $\Om$ if
$Var_H(u;\Om) <\infty$. The space $BV_H(\Om)$ of functions with
bounded $H$-variation in $\Om$, endowed with the norm
\[
||u||_{BV_H(\Om)}\ =\ ||u||_{L^1(\Om)}\ +\ Var_H(u;\Om)\ ,
\]
is a Banach space. A fundamental property of the space $BV_H$ is
the following special case of the compactness Theorem 1.28 proved
in \cite{GN}.

\medskip

\begin{thrm}\label{T:compact}
Let $\Om\subset \bG$ be a (PS) (Poincar\'e-Sobolev) domain. The
embedding
\[
i\ :\ BV_H(\Om)\ \hookrightarrow\ L^q(\Om)
\]
is compact for any $1\leq q < Q/(Q-1)$.
\end{thrm}

\medskip

We now recall a special case of Theorem 1.4 in \cite{CDG}.

\medskip

\begin{thrm}\label{T:iso}
Let $\bG$ be a Carnot group with homogeneous dimension $Q$. There
exists a constant $C(\bG) > 0$, such that for every $g_o\in \bG$,
$0<R<R_o$, one has for every $C^1$ domain $E\subset \overline E
\subset B(g_o,R)$
\[
|E|^{(Q-1)/Q}\ \leq\ C\ P_H(E; B(g_o,R))\ .
\]
\end{thrm}

\medskip

To prove Theorem \ref{T:relisoG0} we need to extend Theorem
\ref{T:iso} from bounded $C^1$ domains to arbitrary sets having
locally finite $H$-perimeter. That such extension be possible is
due in part to the following approximation result for functions in
the space $BV_H$, which is contained in Theorem 1.14 in \cite{GN},
see also \cite{FSS1}.

\medskip

\begin{thrm}\label{T:MSforBV}
Let $\Om\subset \bG$ be open, where $\bG$ is a Carnot group. For
every $u\in BV_H(\Om)$ there exists a sequence $\{u_k\}_{k\in N}$
in $C^\infty(\Om)$ such that
\begin{equation}\label{MS1}
u_k\ \to\ u\ \quad\quad \text{in}\quad L^1(\Om)\quad
\text{as}\quad k \to \infty\ ,
\end{equation}
\begin{equation}\label{MS2}
\underset{k\to\infty}{\lim}\ Var_H(u_k;\Om)\ =\ Var_H(u;\Om)\ .
\end{equation}
\end{thrm}

\medskip

We next introduce the notion of $H$-perimeter.

\medskip

\begin{dfn}\label{D:perimeter}
Let $E\subset \bG$ be a measurable set, $\Om$ be an open set. The
$H$-perimeter of $E$ with respect to $\Om$ is defined by
\[
P_H(E;\Om)\ =\ Var_H(\chi_E;\Om)\ ,
\]
where $\chi_E$ denotes the indicator function of $E$. We say that
$E$ is a $H$-Caccioppoli set if $\chi_E \in BV_H(\Om)$ for every
$\Omega \subset \subset \bG$.
\end{dfn}

\medskip

The reader will notice that when the step of the group $\bG$ is
$r=1$, and therefore $\bG$ is Abelian, the space $BV_H$ coincides
with the space $BV$ introduced by De Giorgi, see \cite{DG1},
\cite{DG2}, \cite{DCP}, and thereby in such setting the Definition
\ref{D:perimeter} coincides with his notion of perimeter. A
fundamental rectifiability theorem \'a la De Giorgi for
$H$-Caccioppoli sets has been established, first for the
Heseinberg group $\Hn$, and then for every Carnot group of step
$r=2$, in the papers \cite{FSS2}, \cite{FSS3}, \cite{FSS4}. We
will need the following simple fact.

\medskip

\begin{lemma}\label{L:equalper}
Let $R_o>0$ be given and consider a $H$-Caccioppoli set $E\subset
\overline E \subset B(e,R_o)$, then
\begin{equation}\label{equalper} P_H(E,B(e,R_o))\ =\
P_H(E,\bG)\ .
\end{equation}
\end{lemma}

\begin{proof}[\textbf{Proof}]
This can be easily seen as follows. Clearly, one has trivially
$P_H(E,B(e,R_o)) \leq P_H(E,\bG)$. To establish the opposite
inequality, let $r_o<R_o$ be such that $E\subset B(e,r_o)$, and
pick $f\in C_0^\infty(B(e,R_o))$ be such that $0\leq f\leq 1$, and
$f\equiv 1$ on $\overline B(e,r_o))$. If $\zeta\in \mathcal
F(\bG)$, then it is clear that $f \zeta \in
\Gamma^1_0(B(e,R_o);H\bG)$, and that $||f
\zeta||_{L^\infty(B(e,R_o))} \leq 1$, i.e., $f \zeta \in \mathcal
F(B(e,R_o))$. We have
\begin{align*}
& \int_{\bG} \chi_{E}\ div_H \zeta\ dg\ =\ \int_{B(e,R_o)}
\chi_{E}\ f\ div_H \zeta\ dg
\\
& =\ \int_{B(e,R_o)} \chi_{E}\ div_H (f \zeta)\ dg\ -\
\int_{B(e,R_o)} \chi_{E}\ <\nabla_H f , \zeta>\ dg
\\
& =\ \int_{B(e,R_o)} \chi_{E}\ div_H (f \zeta)\ dg\ \leq\
P_H(E,B(e,R_o))\ .
\end{align*}

Taking the supremum over all $\zeta \in \mathcal F(\bG;H\bG)$ we
reach the conclusion $P_H(E,B(e,R_o)) \geq P_H(E,\bG)$, thus
obtaining \eqref{equalper}.

\end{proof}

\medskip

In the next result we extend the isoperimetric inequality from
$C^1$ to bounded $H$-Caccioppoli sets.

\medskip

\begin{thrm}\label{T:isoCac}
Let $\bG$ be a Carnot group with homogeneous dimension $Q$. There
exists a constant $C_{iso}(\bG) > 0$ such that for every bounded
$H$-Caccioppoli set $E\subset \bG$ one has
\[
|E|^{(Q-1)/Q}\ \leq\ C_{iso}(\bG)\ P_H(E; \bG)\ .
\]
\end{thrm}

\begin{proof}[\textbf{Proof}]
In \cite{CDG} it was proved that Theorem \ref{T:iso} implies the
following Sobolev inequality of Gagliardo-Nirenberg type: for
every $u\in C^1_0(B(g_o,R))$
\begin{equation}\label{GN}
\left\{\int_{B(g_o,R)}\ |u|^{Q/(Q-1)}\ dg\right\}^{(Q-1)/Q}\ \leq\
C\ \frac{R}{|B(g_o,R)|^{1/Q}}\ \int_{B(g_o,R)}\ |\nabla_H u|\ dg\
.
\end{equation}

If now $u\in BV_H(B(g_o,R))$, with $supp\ u \subset B(g_o,R)$,
then by Theorem \ref{T:MSforBV} there exists a sequence
$\{u_k\}_{k\in \mathbb N}\in C^\infty_0(B(g_o,R))$ such that
$u_k\to u$ in  $L^1(B(g_o,R))$, and $Var_H(u_k;B(g_o,R))\to
Var_H(u;B(g_o,R))$, as $k\to \infty$. Passing to a subsequence, we
can assume that $u_k(g) \to u(g)$, for $dg$-a.e. $g\in B(g_o,R)$.
Applying \eqref{GN} to $u_k$ and passing to the limit we infer
from the theorem of Fatou
\[
\left\{\int_{B(g_o,R)}\ |u|^{Q/(Q-1)}\ dg\right\}^{(Q-1)/Q}\ \leq\
C\ \frac{R}{|B(g_o,R)|^{1/Q}}\ Var_H(u;B(g_o,R))\ ,
\]
for every $u \in BV_H(B(g_o,R))$, with $supp\ u \subset B(g_o,R)$.
If now $E\subset \overline E \subset B(g_o,R)$ is a
$H$-Caccioppoli set, then taking $u = \chi_E$ in the latter
inequality we obtain
\[
|E|^{(Q-1)/Q}\ \leq\ C\ \frac{R}{|B(g_o,R)|^{1/Q}}\ P_H(E;
B(e,R_o))\ .
\]

At this point, to reach the desired conclusion we only need to use
Lemma \ref{L:equalper} and observe that $|B(g_o,R)| = R^Q
|B(e,1)|$. We thus obtain the conclusion with $C_{iso}(\bG) = C
|B(e,1)|^{-1/Q}$.

\end{proof}

\medskip

The following is a basic consequence of Theorem \ref{T:isoCac}.

\medskip

\begin{thrm}\label{T:isoCac2}
Let $\bG$ be a Carnot group with homogeneous dimension $Q$. With
$C_{iso}(\bG)$ as in Theorem \ref{T:iso}, one has for any bounded
$H$-Caccioppoli set
\[
|E|^{(Q-1)/Q}\ \leq\ C_{iso}(\bG)\ P_H(E; \bG)\ .
\]
\end{thrm}

\medskip

To establish Theorem \ref{T:relisoG0} we next prove that one can
remove from Theorem \ref{T:isoCac2}, without altering the constant
$C_{iso}(\bG)$, the restriction that the $H$-Caccioppoli set be
bounded. We recall a useful representation formula. In what
follows $N$ indicates the topological dimension of $\bG$, and
$H_{N-1}$ the $(N-1)$-dimensional Hausdorff measure constructed
with the Riemannian distance of $\bG$.

\medskip

\begin{lemma}\label{L:smoothsets}
Let $\Om \subset \bG$ be an open set and $E\subset \bG$ be a $C^1$
bounded domain. One has
\[
P_H(E;\Om)\ =\ \int_{\Om \cap \partial E} \frac{|\bN_H|}{|\bN|}\
dH_{N-1}\ ,
\]
where $\bN_H = \sum_{j=1}^m <\bN,X_j> X_j$ is the projection onto
$H\bG$ of the Riemannian normal $\bN$ exterior to $E$. In
particular, when
\begin{equation}\label{E1}
E\ =\ \{g\in \bG \mid \phi(g)<0\}\ ,
\end{equation}
with $\phi \in C^1(\bG)$, and $|\nabla \phi| \geq \alpha >0$ in a
neighborhood of $\p E$, then $\bN = \nabla \phi$, and therefore
$|\bN_H| = |\nabla_H \phi|$. When $\Om = \bG$ we thus obtain in
particular \begin{equation}\label{rep}
P_H(E;\bG)\ =\ \int_{\p E}
\frac{|\nabla_H \phi|}{|\nabla \phi|}\ dH_{N-1}\ .
\end{equation}
\end{lemma}

For the proof of this lemma we refer the reader to \cite{CDG}. For
a detailed study of the perimeter measure in Lemma
\ref{L:smoothsets} and \eqref{rep}, we refer the reader to
\cite{DGN1}, \cite{DGN2} and \cite{CG}. We can finally provide the
proof of Theorem \ref{T:relisoG0}.

\medskip

\begin{thrm}\label{T:relisoG}
Let $\bG$ be a Carnot group with homogeneous dimension $Q$. With
the same constant $C_{iso}(\bG)>0$ as in Theorem \ref{T:isoCac2},
for every $H$-Caccioppoli set $E\subset \bG$ one has
\[
|E|^{(Q-1)/Q}\ \leq\ C_{iso}(\bG)\ P_H(E;\bG)\ .
\]
\end{thrm}

\begin{proof}[\textbf{Proof}]
In view of Theorem \ref{T:isoCac2} we only need to consider the case
of an unbounded $H$-Caccioppoli set $E$. If $P_H(E;\bG) = + \infty$
there is nothing to prove, so we assume that $P_H(E;\bG) < + \infty$
and $|E| < + \infty$. We consider the $C^\infty$ $H$-balls $B_H(e,R)
= \{g\in \bG \mid \rho(g)<R\}$, generated by the pseudo-distance
$\rho = \rho_e = \Gamma(\cdot,e)^{1/(2-Q)}\in
C^\infty(\bG\setminus\{e\}) \cap C(\bG)$, where $\Gamma(\cdot,e)\in
C^\infty(\bG\setminus\{e\})$ is the fundamental solution with
singularity at the identity for the sub-Laplacian $\Delta_H =
\sum_{j=1}^m X_j^2$ (the reader should notice that any other smooth
gauge would do). For any $R>0$ we have
\begin{align}\label{proof1}
P_H(E\cap B_H(e,R);\bG)\ & \leq\ P_H(E;B_H(e,R))\ +\
P_H(B_H(e,R);E)\ .
\end{align}

Here, when we write $P_H(B_H(e,R);E)$ we mean the standard measure
theoretic extension of the $H$-perimeter from open sets to Borel
sets, see for instance \cite{Z}. Thanks to the smoothness of
$B_H(e,R)$ we have from Lemma \ref{L:smoothsets}
\[
P_H(B_H(e,R);E)\ =\ \int_{\partial B_H(e,R) \cap E}
\frac{|\bN_H|}{|\bN|}\ dH_{N-1}\ =\ \int_{\partial B_H(e,R) \cap
E} \frac{|\nabla_H \rho|}{|\nabla \rho|}\ dH_{N-1}\ .
\]

Recalling that $\Gamma(\cdot,e)$ is homogeneous of degree $2-Q$,
see \cite{F1}, \cite{F2}, and therefore $\rho$ is homogeneous of
degree one, we infer that for some constant $C(\bG)>0$,
\begin{equation}\label{F3}
|\nabla_H\rho|\ \leq\ C(\bG)\ .
\end{equation}

This gives
\begin{equation}\label{proof2}
P_H(B_H(e,R);E)\ \leq\ C(\bG)\ \int_{\partial B_H(e,R) \cap E}
\frac{dH_{N-1}}{|\nabla \rho|}\ .
\end{equation}

By Federer co-area formula \cite{Fe}, we obtain
\[
\infty\ >\ |E|\ =\ \int_{\bG} \chi_E\ dg\ =\ \int_0^\infty\
\int_{\partial B_H(e,t)\cap E} \frac{dH_{N-1}}{|\nabla \rho|}\ dt\
,
\]
therefore there exists a sequence $R_k \nearrow \infty$ such that
\begin{equation}\label{proof3}
\int_{\partial B_H(e,R_k)\cap E} \frac{dH_{N-1}}{|\nabla \rho|}\
\underset{k \to \infty}{\longrightarrow}\ 0\ .
\end{equation}

Using \eqref{proof3} in \eqref{proof2} we find
\begin{equation}\label{proof4}
\underset{k\to \infty}{\lim}\ P_H(B_H(e,R_k);E)\ =\ 0\ .
\end{equation}

From \eqref{proof1}, \eqref{proof4}, we conclude
\begin{equation}\label{proof5}
\underset{k\to \infty}{\limsup}\ P_H(E\cap B_H(e,R_k);\bG)\ \leq\
P_H(E;\bG)\ .
\end{equation}

We next apply Theorem \ref{T:isoCac2} to the bounded
$H$-Caccioppoli set $E\cap B_H(e,R_k)$ obtaining
\[
|E\cap B_H(e,R_k)|^{(Q-1)/Q}\ \leq\ C_{iso}(\bG)\ P_H(E\cap
B_H(e,R_k);\bG)\ .
\]

Letting $k\to \infty$ in the latter inequality, from
\eqref{proof5}, and from the relation
\[
\underset{k\to \infty}{\lim}\ |E\cap B_H(e,R_k)|^{(Q-1)/Q}\ =\
|E|^{(Q-1)/Q}\ ,
\]
 we conclude that
\[
|E|^{(Q-1)/Q}\ \leq\ C_{iso}(\bG)\ P_H(E;\bG)\ .
\]

This completes the proof.

\end{proof}

\medskip

We close this section with two basic properties of the
$H$-perimeter which clearly play a role also in Theorem
\ref{T:isoine}.

\medskip

\begin{prop}\label{P:blowupP}
In a Carnot group $\bG$ one has for every measurable set $E\subset
\bG$ and every $r>0$
\begin{equation}\label{Pscaling}
P_H(E;\bG)\ =\ r^{Q-1}\ P_H(\delta_{1/r} E ; \bG)\ .
\end{equation}
\end{prop}

\begin{proof}[\textbf{Proof}]
Let $E\subset \bG$ be a measurable set. If $\zeta \in
C^1_0(\bG,H\bG)$, then the divergence theorem, and a rescaling,
give
\begin{equation}\label{ibyparts1}
\int_E\ div_H \zeta\ dg\ =\ \int_E\ \sum_{j=1}^m X_j \zeta_j\ dg\
=\ r^Q\ \int_{E_r}\ \sum_{j=1}^m X_j \zeta_j(\delta_r g)\ dg\ ,
\end{equation}
where we have let $E_r = \delta_{1/r}(E) = \{g\in \bG \mid
\delta_r g \in E \}$. Since
\[
X_j(\zeta_j \circ \delta_r)\ =\ r\ (X_j \zeta_j \circ \delta_r)\ ,
\]
we conclude
\begin{equation}\label{ibyparts2}
\int_E\ \sum_{j=1}^m X_j \zeta_j\ dg\ =\ r^{Q-1}\ \int_E\
\sum_{j=1}^m X_j(\zeta_j\circ \delta_r)\ dg\ .
\end{equation}

Formula \eqref{ibyparts2} implies the conclusion.

\end{proof}

\medskip

Proposition \ref{P:blowupP} asserts that the $H$-perimeter scales
appropriately with respect to the non-isotropic group dilations.
Since on the other hand one has $|\delta_{1/r} E| = r^{-Q} |E|$,
we easily obtain the following important scale invariance of the
isoperimetric quotient.

\medskip

\begin{prop}\label{P:invariance}
For any $H$-Caccioppoli set in a Carnot group $\bG$ one has
\begin{equation}
\frac{P_H(E;\bG)}{|E|^{(Q-1)/Q}}\ =\ \frac{P_H(\delta_{1/r}
E;\bG)}{|\delta_{1/r} E|^{(Q-1)/Q}}\ ,\quad\quad\quad r > 0\ .
\end{equation}
\end{prop}

\medskip

Another equally important fact, which is however a trivial
consequence of the left-invariance on the vector fields
$X_1,...,X_m$, and of the definition of $H$-perimeter, is the
translation invariance of the isoperimetric quotient.

\medskip

\begin{prop}\label{P:invariance2}
For any $H$-Caccioppoli set in a Carnot group $\bG$ one has
\begin{equation}
\frac{P_H(L_{g_o}(E) ;\bG)}{|L_{g_o} (E)|^{(Q-1)/Q}}\ =\
\frac{P_H(E;\bG)}{| E|^{(Q-1)/Q}}\ ,\quad\quad\quad g_o \in \bG\ ,
\end{equation}
where $L_{g_o}g = g_o g$ is the left-translation on the group.
\end{prop}

\vskip 0.6in


\section{\textbf{Partial solution of the isoperimetric problem in $\Hn$ }}\label{S:PS}

\vskip 0.2in

The objective of this section is proving Theorems
\ref{T:isoprofile} and \ref{T:isoine}. This will be accomplished
in several steps. First, we introduce the relevant notions and
establish some geometric properties of the $H$-perimeter that are
relevant to the isoperimetric profiles. Next, we collect some
results from convex analysis and calculus of variations. Finally,
we proceed to proving Theorems \ref{T:isoprofile} and
\ref{T:isoine}. In what follows we adopt the classical
non-parametric point of view, see for instance \cite{MM},
according to which a $C^2$ hypersurface $\mS \subset \bG$ locally
coincides with the zero set of a real function. Thus, for every
$g_0\in \mS$ there exists an open set $\mathcal O\subset \bG$ and
a function $\phi\in C^2(\mathcal O)$ such that: (i) $|\nabla
\phi(g)| \not= 0$ for every $g\in \mathcal O$; (ii) $\mS\cap
\mathcal O = \{g\in \mathcal O\mid \phi(g) = 0\}$. We will always
assume that $\mS$ is oriented in such a way that for every $g\in
\mS$ one has
\[
\bN(g) = \nabla \phi(g) = X_1 \phi(g) X_1 + ...  + X_m \phi(g) X_m
+ ... + X_{r,1} \phi(g) X_{r,1} + ... + X_{r,m_r}\phi(g)
X_{r,m_r}\ .
\]

To justify the second equality the reader should bear in mind that
we have endowed $\bG$ with a left-invariant Riemannian metric with
respect to which $\{X_1,...,X_m,...,X_{r,m_r}\}$ constitute an
orthonormal basis. Given a surface $\mS \subset \bG$, we let
\begin{equation}\label{pi}
 p_i\ =\ <\bN,X_i>\ , \quad\quad\quad i =
1,...,m\ , \end{equation}
 and define the angle function
\begin{equation}\label{W}
 W\ =\ \sqrt{p_1^2 + ... + p_m^2}\ . \end{equation}

The motivation for the name comes from the fact that, if $U\angle
V$ denotes the angle between two vector fields $U,V$ on $\bG$,
then
\begin{equation}\label{angle}
\cos(\nuX\angle \bN)\ =\ \frac{<\nuX,\bN>}{|\bN|}\ =\
\frac{W}{|\bN|}\ .
\end{equation}

The characteristic locus of $\mS$ is the closed set
\[
\Sigma\ =\ \{g\in \mS\mid W(g) = 0\}\  =\ \{g\in \mS \mid H_g\bG
\subset T_g \mS\}\ . \]

We recall that is was proved in \cite{B}, \cite{Ma} that $\mathcal
H^{Q-1}(\Sigma) = 0$, where $\mathcal H^s$ denotes the
$s$-dimensional Hausdorff measure associated with the
Carnot-Carath\'eodory distance of $\bG$, and $Q$ indicates the
homogeneous dimension of $\bG$. We also recall the earlier result of
Derridj \cite{De1}, \cite{De2}, which states that when $\mS$ is
$C^\infty$ the standard surface measure of $\Sigma$ vanishes. Later
on in this section we will need a result from \cite{B}, see Theorem
\ref{T:balogh} below.

On the set $\mathcal S\setminus \Sigma$ we define the
\emph{horizontal Gauss map} by
\begin{equation}\label{gm}
\nuX\ =\  \pb_1 X_1 + ... + \pb_m X_m\ ,
\end{equation}
where we have let
\begin{equation}\label{pbar}
\pb_1\ =\ \frac{p_1}{W}\ ,\ ...\ ,\ \pb_m\ =\ \frac{p_m}{W}\ ,
\quad\quad\text{so that}\quad\quad |\nuX|^2\ =\ \pb_1^2 + ... +
\pb_m^2\ \equiv\ 1\quad\quad \text{on}\quad\quad \mathcal S
\setminus \Sigma\ .
\end{equation}

Given a point $g_0\in \mS\setminus \Sigma$, the horizontal tangent
space of $\mS$ at $g_0$ is defined by \[ T_{H,g_0}(\mS)\ =\
\{\boldsymbol v \in H_{g_0}\bG \mid <\boldsymbol v,\nuX(g_0)>\ =\
0\}\ . \]

For instance, when $\bG = \HH$, then a basis for $T_{H,g_0}(\mS)$
is given by the single vector field
\begin{equation}\label{nup} \nup\ =\ \pb_2\ X_1\ -\
\pb_1\ X_2\ .
\end{equation}

Given a function $u\in C^1(\mS)$ one clearly has $\delta_H
u(g_0)\in T_{H,g_0}(\mS)$. We next recall some basic definitions
from \cite{DGN3}.

\medskip

Let $\nabla^H$ denote the horizontal Levi-Civita connection
introduced in \cite{DGN3}. Let $\mS\subset \bG$ be a $C^2$
hypersurface. Inspired by the Riemannian situation we introduce a
notion of horizontal second fundamental on $\mS$ as follows.

\medskip

\begin{dfn}\label{D:sff}
Let $\mS\subset \bG$ be a $C^2$ hypersurface, with $\Sigma =
\varnothing$, then we define a tensor field of type $(0,2)$ on
$T_H\mS$, as follows: for every $X,Y\in C^1(\mS;T_H\mS)$
\begin{equation}\label{sff} II^{H,\mS}(X,Y)\ =\ <\nabla^H_X Y,\nuX>
\nuX\ .
\end{equation}
We call $II^{H,\mS}(\cdot,\cdot)$ the \emph{horizontal second
fundamental form} of $\mS$. We also define $\mathcal A^{H,\mS} : T_H
\mS \to T_H \mS$ by letting for every $g\in \mS$ and $\boldsymbol u,
\boldsymbol v \in T_{H,g}$
\begin{equation}\label{shape}
<\mathcal A^{H,\mS} \boldsymbol u,\boldsymbol v>\ =\ -\
<II^{H,\mS}(\boldsymbol u,\boldsymbol v),\nuX>\ =\ -\ <\nabla_X^H
Y,\nuX>\ ,
\end{equation}
where $X, Y \in C^1(\mS,T_H\mS)$ are such that $X_g = \boldsymbol
u$, $Y_g = \boldsymbol v$. We call the endomorphism $\mathcal
A^{H,\mS} : T_{H,g}\mS \to T_{H,g}\mS$ the \emph{horizontal shape
operator}. If $\boldsymbol e_1,...,\boldsymbol e_{m-1}$ denotes a
local orthonormal frame for $T_H\mS$, then the matrix of the
horizontal shape operator with respect to the basis $\boldsymbol
e_1,...,\boldsymbol e_{m-1}$ is given by the $(m-1)\times(m-1)$
matrix $ \big[- <\nabla_{\boldsymbol e_i}^H \boldsymbol
e_j,\nuX>\big]_{i,j=1,...,m-1}$.

\end{dfn}

\medskip

By the horizontal Koszul identity in \cite{DGN3}, one easily
verifies that
\[
<\nabla_{\boldsymbol e_i}^H \boldsymbol e_j,\nuX>\ =\ -\
<\nabla_{\boldsymbol e_i}^H \nuX, \boldsymbol e_j>\ .
\]

Using Definition \ref{D:sff} one recognizes that
\begin{equation}\label{nonsymm}
II^{H,\mS}(X,Y)\ -\ II^{H,\mS}(Y,X)\ =\ <[X,Y]^H,\nuX> \nuX\ ,
\end{equation}
and therefore, unlike its Riemannian counterpart, the horizontal
second fundamental form of $\mS$ is not necessarily symmetric. This
depends on the fact that, if $X,Y\in C^1(\mS;HT\mS)$, then it is not
necessarily true that the projection of $[X,Y]$ onto the horizontal
bundle $H\Hn$, $[X,Y]^H$, belongs to $C^1(\mS;T_H\mS)$.

\medskip

\begin{dfn}\label{D:HMC}
We define the \emph{horizontal principal curvatures} as the real
eigenvalues $\kappa_1,...,\kappa_{m-1}$ of the symmetrized operator
\[
\mathcal A^{H,\mS}_{sym}\ =\ \frac{1}{2}\left\{\mathcal A^{H,\mS} +
(\mathcal A^{H,\mS})^t\right\}\ ,
\]
The $H$-mean curvature of $\mS$ at a non-characteristic point
$g_0\in \mS$ is defined as
\[
\mathcal H\ =\ -\ trace\ \mathcal A^{H,\mS}_{sym}\ =\
\sum_{i=1}^{m-1} \kappa_i\ =\ \sum_{i=1}^{m-1}
<\nabla^H_{\boldsymbol e_i} \boldsymbol e_i,\nuX>\ .
\]
If $g_0$ is characteristic, then we let
\[
\mathcal H(g_0)\ =\ \underset{g\to g_0, g\in \mathcal S\setminus
\Sigma}{\lim}\ \mathcal H(g)\ ,
\]
provided that such limit exists, finite or infinite. We do not
define the $H$-mean curvature at those points $g_0\in \Sigma$ at
which the limit does not exist. Finally, we call $\vec{\mathcal H} =
\mathcal H \nuX$ the $H$-\emph{mean curvature vector}.
\end{dfn}

\medskip

Hereafter, when we say that a function $u$ belongs to the class
$C^k(\mS)$, we mean that $u\in C(\mS)$ and that for every $g_0\in
\mS$, there exist an open set $\mathcal O\subset \HH$, such that $u$
coincides with the restriction to $\mS\cap \mathcal O$ of a function
in $C^k(\mathcal O)$. The tangential horizontal gradient of a
function $u\in C^1(\mS)$ is defined as follows
 \begin{equation}\label{del}
\nabla^{H,\mS} u\ =\ \nabh u\ -\ <\nabh u,\nuX> \nuX\ .
\end{equation}

The definition of $\nabla^{H,\mS} u$ is well-posed since it is noted
in \cite{DGN3} that it only depends on the values of $u$ on $\mS$.
Since $|\nuX| \equiv 1$ on $\mS \setminus \Sigma$, we clearly have
$<\nabla^{H,\mS} u,\nuX> = 0$, and therefore
\begin{equation}\label{pitagora}
|\nabla^{H,\mS} u|^2\ =\ |\nabh u|^2\ -\ <\nabh u,\nuX>^2\ .
\end{equation}

\medskip

\begin{dfn}\label{D:cmc}
We say that a $C^2$ hypersurface $\mathcal S$ has \emph{constant
$H$-mean curvature} if $\mathcal H$ is globally defined on
$\mathcal S$, and $\mathcal H \equiv const.$ We say that $\mathcal
S$ is $H$-\emph{minimal} if $\mathcal H \equiv 0$.
\end{dfn}

\medskip

Minimal surfaces have been recently studied in \cite{Pa}, \cite{GP},
\cite{CHMY}, \cite{CH}, \cite{DGN5}, \cite{DGNP}, \cite{BSV}. The
last two papers contain also a complete solution of the Bernstein
type problem for the Heisenberg group $\HH$. The following result is
taken from \cite{DGN3}.

\medskip

\begin{prop}\label{P:equalMC}
The $H$-mean curvature coincides with the function
\begin{equation}\label{equal2}
\mathcal H\ =\ \sum_{i=1}^m\ \nabla^{H,\mS} \pb_i\ =\ \sum_{i=1}^m\
X_i \pb_i\ .
\end{equation}
\end{prop}

For instance, when $\bG = \HH$, then according to Proposition
\ref{P:equalMC}, the $H$-mean curvature of $\mathcal S$ is given
by
\begin{equation}\label{HH}
\mathcal H\ =\ \sum_{i=1}^2 \nabla^{H,\mS}_i \nui\ =\
\nabla^{H,\mS}_1(\pb_1)\ +\ \nabla^{H,\mS}_2(\pb_2)\ =\ X_1 \pb_1\
+\ X_2 \pb_2\ ,\quad\quad\quad\quad\text{on}\quad \mS \setminus
\Sigma\ .
\end{equation}

In this situation, given a $C^2$ surface $\mS \subset \HH$, there is
only one horizontal principal curvature $\kappa_1(g_0)$ at every
$g_0\in \mS \setminus \Sigma$. Since in view of \eqref{nup} the
vector $\nup(g_0)$ constitutes an orthonormal basis of
$T_{H,g_0}(\mS)$, according to Definition \ref{D:sff} we have
\[ \kappa_1(g_0)\ =\ II_H(\nup(g_0),\nup(g_0))\ .
\]

One can verify, see \cite{DGN3}, that the right-hand side of the
latter equation equals $- \mathcal H(g_0)$.  We recall one more result
concerning the $H$-mean curvature that will be useful in the proof of
Proposition \ref{P:necessary}.  Details can be found in \cite{DGN3}.

\begin{prop}\label{P:tphi}
Suppose $\mathcal S \subset \Hn$ is a level set of a function
$\phi$ that takes the form
\[
\phi(z,t)\ =\ t - u\big(\frac{|z|^2}{4}\big)\ ,
\]
for some $C^2$ function $u:[0,\infty) \to \mathbb R$.
For every point point $g = (z,t)\in \mathcal S$ such that $z\not= 0$ the $H$-mean curvature at $g$ is given by
\begin{equation}\label{tphi}
\mathcal H\ =\ -\,\frac{2\ s\ u''(s)\ +\ (Q -3)\ u'(s)\ (1 +
u'(s)^2)}{2\ \sqrt s\ (1 + u'(s)^2)^{3/2}}\ ,\quad\quad s =
\frac{|z|^2}{4}\ .
\end{equation}
\end{prop}

\medskip

In Proposition \ref{P:tphi} the hypothesis $z\not= 0$ is justified
the fact that, under the given assumptions, if $\mS$ intersects
the $t$-axis in $\Hn$, then the points of intersections are
necessarily characteristic for $\mS$.

Hereafter in this paper, we restrict our attention to $\bG = \Hn$.
In Definition \ref{D:cmc}, following the classical tradition, we
have called a hypersurface $H$-minimal if its $H$-mean curvature
vanishes identically. However, in the classical setting the
measure theoretic definition of minimality is also based on the
notion of local minimizer of the area functional. In the paper
\cite{DGN3} we have proved that there is a corresponding
sub-Riemannian counterpart of such interpretation based on
appropriate first and second variation formulas for the
$H$-perimeter. For instance, the following first variation formula
holds in the Heisenberg group $\mathbb H^1$.

\medskip

 \begin{thrm}\label{T:variations}
Let $\mS\subset \HH$ be an oriented $C^2$ surface, then the first
variation of the $H$-perimeter with respect to the deformation
\begin{equation}\label{def}
J_\lambda(g) = g + \lambda \mathcal X(g) = g + \lambda \big(a(g)
X_1 + b(g) X_2 + k(g) T\big)\ ,\quad\quad\quad g = (x,y,t)\in \mS\
,
\end{equation}
is given by
\begin{equation}\label{fvH}
\frac{d}{d\lambda} P_H(\mathcal S^\lambda)\Bigl|_{\lambda = 0}\ =\
 \int_{\mathcal S}
\mathcal H\ \frac{\cos(\mathcal X \angle \bN)}{\cos(\nuX \angle
\bN)}\ |\mathcal X|\ d\sigma_H\ ,
\end{equation}
where $\angle$ denotes the angle between vectors in the inner
product $<\cdot,\cdot>$. In particular, $\mathcal S$ is stationary
with respect to \eqref{def} if and only if it is $H$-minimal.
\end{thrm}

\medskip

Versions of Theorem \ref{T:variations} have also been obtained
independently by other people. An approach based on motion by
$H$-mean curvature can be found in \cite{BC}. When $a = \pb h$, $b
= \qb h$, and $h\in C^\infty_0(\mS\setminus \Sigma)$, then a proof
based on CR-geometry can be found in \cite{CHMY}. A Riemannian
geometry proof, valid in any $\Hn$, can be found in \cite{RR}.

In what follows we set
\[
\Hn_+\ =\ \{(z,t)\in \Hn\mid t>0\}\ ,\quad\quad \Hn_{-}\ =\
\{(z,t)\in \Hn\mid t<0\}\ . \]

 Consider a domain $\Om \subset
\mathbb R^{2n}$ and a $C^1$ function $u:\Om \to [0,\infty)$. We
assume that $E\subset \Hn$ is a $C^1$ domain for which
\[
E\cap \Hn_+\ =\ \{(z,t)\in \Hn\,|\, z\in \Om,\, 0 < t < u(z)\}\ .
\]

The reader should notice that, since $u>0$ in $\Om$, the graph of
$u$ is not allowed to have flat parts. For $z = (x,y)\in \mathbb
R^{2n}$, we set $z^\perp = (y,-x)$. Indicating with $\phi(z,t) = t -
u(z)$ the defining function of $E\cap \Hn_+$, a simple computation
gives
\begin{align}\label{per1}
|\nH\phi|=\ & =\ \sqrt{\left|\nabla_x u + \frac{y}{2}\right|^2 +
\left|\nabla_y u - \frac{x}{2}\right|^2}
\\
& =\ \left|\nabla_z u + \frac{z^\perp}{2}\right|\ . \notag
\end{align}

The reader should be aware that in the latter equation, the norm in
the left-hand side comes from the Riemannian inner product on $T\Hn
\cong \Hn$, whereas the norm in the right-hand side is simply the
Euclidean norm in $\R^{2n}$. Invoking the representation formula
\eqref{rep} for the $H$-perimeter, which presently gives
\[
P_H(E;\Hn_+)\ =\ \int_{\p E \cap \Hn_+}\ \frac{|\nH\phi|}{|\nabla
\phi|}\ dH_{2n}\ ,
\]
and keeping in mind that, see \eqref{per1}, $|\nabla \phi| = \sqrt{1
+ |\nH\phi|^2}$, and that $dH_{2n} = \sqrt{1 + |\nH \phi|^2}\ dz$,
we obtain
\begin{equation}\label{E:X-per}
P_H(E;\Hn_+)\ =\ \int_\Om \left|\nabla_z u +
\frac{z^\perp}{2}\right|\ dz\ =\ \int_\Om \sqrt{|\nabla_z u|^2 +
\frac{|z|^2}{4} + <\nabla_z u,z^\perp>}\ \ dz\ .
\end{equation}

When $F\subset \Hn$ is a closed set we define
\[
P_H(E;F)\ = \ \underset{F\subset \Om ,   \Om\ \text{open}}{\inf}\
P_H(E;\Om) \ .
\]

Let now $u\in C^1(\Om)$, $u\geq 0$, then using the latter formula we
obtain the following generalization of \eqref{E:X-per}
\begin{equation}\label{perclosed}
P_H(E;\overline{\Hn_+})\ =\ \int_\Om \left|\nabla_z u +
\frac{z^\perp}{2}\right|\ dz\ =\ \int_\Om \sqrt{|\nabla_z u|^2 +
\frac{|z|^2}{4} + <\nabla_z u,z^\perp>}\ \ dz\ .
\end{equation}

The reader should notice that, unlike \eqref{E:X-per}, in equation
\eqref{perclosed} we allow the graph of $u$ to have flat parts,
i.e., subsets of $\Om$ in which the function $u$ vanishes.

 In what follows, we
recall an invariance property of the $H$-perimeter which plays a
role in the proof of Theorem \ref{T:isoprofile}. Consider the map
$\mathcal O: \Hn \to \Hn$ defined by
\[
\mathcal O(x,y,t)\ =\ (y,x,-t)\ .
\]

It is obvious that $\mathcal O$ preserves Lebesgue measure (which is
a bi-invariant Haar measure on $\Hn$). In fact, the map $\mathcal O$
is a group and Lie algebra automorphism of $\Hn$. Such map is called
\emph{inversion} in \cite{F3}, p.20. Using the properties of the map
$\mathcal O$ and a standard contradiction argument, one can easily
prove the following result.

\medskip

\begin{thrm}\label{T:iso-symm}
Let $E\subset \Hn$ be a bounded open set such that $\p E \cap
\Hn_+$ and $\p E \cap \Hn_-$ are $C^1$ hypersurfaces, and assume
that $E$ satisfies the following condition: there exists $R
> 0$ such that
\begin{equation}\label{E:red-partial-symm}
E\ \cap\ \{t=0\}\ =\  B(0,R) \ .
\end{equation}
Suppose $E$ is an isoperimetric set satisfying $|E\cap \Hn_+|\ =\
|E\cap \Hn_-|\ =\ |E|/2$, then
\[
P_H(E;\overline{\Hn_+})\ =\ P_H(E;\overline{\Hn_-})\ .
\]
\end{thrm}

\medskip

We now introduce the relevant functional class for our problem. The
space of competing functions $\mathcal D$ is defined as follows.
Consider the vector space $\mathcal V = C_0(\R^{2n})$.

\medskip

\begin{dfn}\label{D:fc}
We let
\begin{align}\label{E:D}
\mathcal D = & \{u\in \mathcal V |\text{there exists}\ R > 0\
\text{such that}\ u\geq 0 \text{ in } B(0,R),
\\
& \overline{B}(0,R) = \bigcap \{B(0,R+ \rho)\mid\ supp(u)\subset
B(0,R + \rho)\}\notag\\ & u\in C^{1,1}_{loc}(B(0,R))\cap
W^{1,1}(B(0,R))\}\ . \notag
\end{align}
\end{dfn}

We note explicitly that, as a consequence of Definition
\eqref{D:fc}, if $u\in \mathcal D$ and $R$ is as in \eqref{E:D}, we
have $u=0$ on $\p B(0,R)$. Furthermore, the functions in $\mathcal
D$ are allowed to have large sets of zeros, i.e., their graph is
allowed to touch the hyperplane $t=0$ in sets of large measure. We
remark that $\mathcal D$ is not a vector space, nor it is a convex
subset of $\mathcal V$. We mention that the requirement $u\in
C^{1,1}_{loc}(B(0,R))$ in the definition of the class $\mathcal D$,
is justified by the following considerations. When we compute the
Euler-Lagrange equation of the $H$-perimeter functional
\eqref{E:X-per} we need to know that, with $\Om = B(0,R)$, the set
$\{z = (x,y)\in  \Om\subset \R^{2n}\mid |\nabla_z u(z) +
\frac{z^\perp}{2}| = 0\}$, which is the projection of the
characteristic set of the graph of $u$ onto $\R^{2n}\times \{0\}$,
has vanishing $2n$-dimensional Lebesgue measure. This is guaranteed
by the following sharp result of Z. Balogh (see Theorem 3.1 in
\cite{B}) provided that $u\in C^{1,1}_{loc}(\Om)$, but it fails in
general for $u\in C^{1,\alpha}_{loc}(\Om)$ for every $0<\alpha<1$.

\medskip

\begin{thrm}\label{T:balogh}
Let $\Om = B(0,R) \subset \R^{2n}$ and consider $u\in
C^{1,1}_{loc}(\Om)$, then $|\mathcal A(u)| = 0$, where $\mathcal
A(u)= \{z \in \Om \mid \nabla_z u(z) + z^\perp/2 = 0\}$, and $|E|$
denotes the $2n$-dimensional Lebesgue measure of $E$ in $\R^{2n}$.
If instead $u\in C^2(\Om)$, then the Euclidean dimension of $E$ is
$\leq n$.
\end{thrm}

\medskip

Following classical ideas from the Calculus of Variation, we next
introduce the admissible variations for the problem at hand, see
\cite{GH} and \cite{Tr}.

\medskip

\begin{dfn}\label{D:admissible}
Given $u\in \mathcal D$, we say that $\phi\in \mathcal V$, with
$supp\ \phi \subseteq supp\ u$, is $\mathcal D$-admissible at $u$ if
$u + \lambda \phi \in \mathcal D$ for all $\lambda\in \R$
sufficiently small .
\end{dfn}

\medskip

Now, for $u\in\mathcal D$ we let
\begin{equation}\label{E:constraint2}
\mathcal G[u]\ =\ \int_{supp(u)} u(z)\, dz\ \ =\ \int_{B(0,R)}
u(z)\, dz\ .
\end{equation}
With \eqref{E:X-per} in mind, we define for such $u$

\begin{equation}\label{E:minimize2}
\mathcal F[u] \ =\ \int_{supp(u)} \sqrt{|\nabla_z u|^2 +
\frac{|z|^2}{4}\ + <\nabla_z u,z^\perp>}\ \ .
\end{equation}

In the class of $C^1$ graphs over $\R^{2n}\times \{0\}$, the
isoperimetric problem consists in minimizing the functional
$\mathcal F[u]$, subject to the constraint that $\mathcal G[u] = V$,
where $V>0$ is given and $B(0,R)$ is replaced by an a priori unknown
domain $\Om$. We emphasize that finding the section of the
isoperimetric profile with the hyperplane $\{t=0\}$, i.e., finding
the domain $\Om$, constitutes here part of the problem. Because of
the lack of an obvious symmetrization procedure, this seems a
difficult question at the moment. To make further progress we
restrict the class of domains $E$ by imposing that their section
with the hyperplane $\{t=0\}$ be a ball, i.e., we assume that, given
$E\in \mathcal E$, there exists $R = R(E)
>0$ such that $\Om = B(0,R)$. Under this hypothesis, we can appeal
to Theorem \ref{T:iso-symm}. The latter implies that it suffices to
solve the following variational problem: \emph{given $V>0$, find
$R_o>0$ and $u_o\in\mathcal D$ with $supp(u_o) = B(0,R_o)$
 for which the following holds}

\begin{equation}\label{E:Gen_Var}
\mathcal F[u_o]\ =\  min\{\mathcal F[u]\,|\, u\in\mathcal D\}
\quad\quad \text{and} \quad\quad \mathcal G[u_o]\ =\ \frac{V}{2}\ .
\end{equation}

To reduce the problem \eqref{E:Gen_Var} to one without constraint,
we will apply the following standard version of the Lagrange
multiplier theorem (see, e.g., Proposition 2.3 in \cite{Tr}).

\medskip

\begin{prop}\label{P:LagrangeM}
Let $\mathcal D$ be a subset of a normed vector space $\mathcal
V$, and consider functionals $\mathcal F$, $\mathcal G_1,\mathcal
G_2$,...,$\mathcal G_k$ defined on $\mathcal D$.  Suppose there
exist constants $\lambda_1,...,\lambda_k \in \R$, and
$u_o\in\mathcal D$, such that $u_o$ minimizes (uniquely)
\begin{equation}\label{E:LagrangeM}
\mathcal F\ +\ \lambda_1\,\mathcal G_1 \ +\ \lambda_2\,\mathcal
G_2 \ +\ \cdots \ +\ \lambda_k\,\mathcal G_k\
\end{equation}
on $\mathcal D$, then $u_o$ minimizes $\mathcal F$ (uniquely) when
restricted to the set
\[
\{u\in\mathcal D\,|\, \mathcal G_j[u]\ =\ \mathcal G_j[u_o], \,
j=1,...,k\}\ .
\]
\end{prop}

\medskip

\begin{rmrk}\label{R:LagrangeM}
The procedure of applying the above proposition to solving a
problem of the type
\[
\begin{cases}
\text{minimize}\quad\{\mathcal F[u]\,|\,u\in\mathcal D\}\ ,
\\
\text{subject to the constraints}\quad\quad\quad \mathcal G_1[u]\
=\ V_1\ , \ \dots\ \ ,\ \mathcal G_k[u]\ =\ V_k\ , \end{cases}
\]
consists of two main steps.  First, one needs to show that
constants $\lambda_1,...,\lambda_k$ and a $u_o\in\mathcal D$ can
be found in such a way that $u_o$ solves the Euler-Lagrange
equation of \eqref{E:LagrangeM}, and $u_o$ satisfies $\mathcal
G_1[u_o] = V_1$,...,$\mathcal G_k[u_o] = V_k$. Finally, one proves
that the solution $u_o$ of the Euler-Lagrange equation is indeed a
minimizer of \eqref{E:LagrangeM}. If the functional involved
possesses appropriate convexity properties, then one can show that
such minimizer $u_o$ is unique.
\end{rmrk}

\medskip

We then proceed with the first step outlined in the Remark
\ref{R:LagrangeM}. In what follows, with $z\in \R^{2n}$, $u\in
\R$, and $p = (p_1,p_2)\in \mathbb R^{2n}$, we let
\begin{equation}\label{E:integrants}
\begin{cases}
f(z,u,p)\ & =\ f(z,p)\ =\ \sqrt{\left|p_1 + \frac{y}{2}\right|^2 +
\left|p_2 - \frac{x}{2}\right|^2}\ =\ \left|p +
\frac{z^\perp}{2}\right|\ ,
\\
g(z,u,p)\ & =\ g(u)\ =\ u\ ,
\\
h(z,u,p)\ & =\ f(z,p)\ +\ \lambda\ g(u)\ .
\end{cases}
\end{equation}

The constrained variational problem \eqref{E:Gen_Var} is then
equivalent to the following one without constraint (provided the
parameter $\lambda$ is appropriately chosen): \emph{to minimize
the functional
\begin{equation}\label{E:unconstrained2}
\mathcal F[u]\ =\ \int_{supp(u)} h(z,u(z),\nabla_z u(z))\,dz\ =\
\int_{supp(u)} \left\{\left|\nabla_z u(z) +
\frac{z^\perp}{2}\right|\ +\ \lambda u(z)\right\}\ dz\ ,
\end{equation}
over the set $\mathcal D$ defined in \eqref{E:D}}. We easily
recognize that the Euler-Lagrange equation of
\eqref{E:unconstrained2} is
\begin{equation}\label{E:EL-2}
div_z \left[\frac{\nabla_z u + \frac{z^\perp}{2}}{\sqrt{|\nabla_z
u|^2 + \frac{|z|^2}{4} + <\nabla_z u,z^\perp>}}\right]\ = \
\lambda\ .
\end{equation}

\medskip

\begin{rmrk}\label{R:interpretation}
Before proceeding we note explicitly that, if $u\in C^2(\Om)$, and
we consider the $C^2$ hypersurface $\mS = \{(z,t)\in \Hn\mid z\in
\Om\ ,\ t = u(z)\}$, indicating with $\Sigma$ its characteristic
set, then $g = (z,t) \not\in \Sigma$ if and only if $|\nabla_z u|^2
+ \frac{|z|^2}{4} + <\nabla_z u,z^\perp> \not= 0$. In this
situation, using Proposition \ref{P:equalMC}, it can be recognized
that, at every $g\not\in \Sigma$, the quantity in the left-hand side
of \eqref{E:EL-2} represents the $H$-mean curvature $\mathcal H$ of
$\mS$.
\end{rmrk}

\medskip

As we have said, solving \eqref{E:EL-2} on an arbitrary domain of
$\Om \subset \R^{2n}$ is a difficult task. However, when $\Om$ is a
ball in $\mathbb R^{2n}$, the equation \eqref{E:EL-2} admits a
remarkable family of spherically symmetric solutions. We note
explicitly that for a graph $t =u(z)$ with spherical symmetry in
$z$, the only characteristic points can occur at the intersection of
the graph with the $t$-axis.

\medskip

\begin{thrm}\label{T:ss0}
Given $R>0$, for every
\begin{equation}\label{l}
-\ \frac{Q-2}{R}\ \leq\ \lambda <\ 0\ ,
\end{equation}
 the equation
\eqref{E:EL-2}, with the Dirichlet condition $u = 0$ on $\p \Om$,
where  $\Om = B(0,R) = \{z\in \mathbb R^{2n} \mid |z|<R\}$, admits
the spherically symmetric solution $u_{R,\lambda} \in \mathcal D\cap
C^2(\Om\setminus \{0\})$, with
\begin{equation}\label{uo0}
u_{R,\lambda}(z)\ =\ C_{R,\lambda}\ +\ \frac{|z|}{4\lambda}\
\sqrt{(Q-2)^2 - (\lambda |z|)^2}\ \ -\
\frac{(Q-2)^2}{4\lambda^2}\,\sin^{-1}\left(\frac{\lambda
|z|}{Q-2}\right)\ ,
\end{equation}
and
\begin{equation}\label{c}
C_{R,\lambda}\ =\ -\ \frac{R}{4\lambda} \sqrt{(Q-2)^2 - (\lambda
R)^2}\ -\ \frac{(Q -2)^2}{4 \lambda^2}\ \sin^{-1}\
\left(\frac{\lambda R}{Q-2}\right)\ .
\end{equation}
\end{thrm}

\begin{proof}[\textbf{Proof}]
We look for a spherically symmetric solution in the form $u(z) =
\rad(|z|^2/4)$, for some function $\rad\in C^2((0,R^2/4))\cap
C([0,R^2/4])$, with $\rad(R^2/4) = 0$. The equation \eqref{E:EL-2}
becomes
\begin{equation}\label{iso2_10}
div_z \left[\frac{\rad'(|z|^2/4)\ z\ +\ z^\perp}{|z|\ \sqrt{1 +
\rad'(|z|^2/4)}}\right]\ = \ \lambda\ ,\quad\quad\text{in}\
B(0,R)\setminus\{0\}\ .
\end{equation}

Since
\[
div_z \left[\frac{z^\perp}{|z|\ \sqrt{1 + \rad'(|z|^2/4)}}\right]\
= \ 0\ ,
\]
we obtain that \eqref{iso2_10} reduces to the equation
\begin{equation}\label{iso3_10}
div_z \left[\frac{\rad'(|z|^2/4)\ z}{|z|\ \sqrt{1 +
\rad'(|z|^2/4)}}\right]\ = \ \lambda\ .
\end{equation}

The transformation
\begin{equation}\label{F0}
F(r)\ \overset{def}{=}\ \frac{\rad
'\left(\frac{r^2}{4}\right)}{r\sqrt{1+\left(\rad
'\left(\frac{r^2}{4}\right)\right)^2}}\ ,
\end{equation}
turns the nonlinear equation \eqref{iso3_10} into the following
linear one
\begin{equation*}\label{s1}
F'(r)\ +\ \frac{2n}{r}\ F(r)\ =\ \frac{\lambda}{r}\ ,
\end{equation*}
which is equivalent to \[ (r^{2n} F)'\ =\ \lambda\ r^{2n-1}\ .
\]

We note that \[ |r^{2n} F(r)|\ \leq\ r^{2n-1}\ , \quad\quad\quad
\text{for}\quad 0<r\leq \frac{R^2}{4}\ , \] therefore we conclude
that $\underset{r\to 0}\lim\ {r^{2n}F(r)} = 0$. We can thus easily
integrate the above ode, obtaining $F(r) \equiv \lambda/2n$. Setting
$s = r^2/4$ in the latter identity one obtains from \eqref{F0}
\begin{equation}\label{iso4_40}
\frac{\rad '(s)}{\sqrt{1+(\rad '(s))^2}}\ =\ \frac{\lambda}{n}\
\sqrt s\ =\ \frac{2 \lambda}{Q-2} \sqrt s\ .
\end{equation}

Excluding the case of $H$-minimal surfaces (corresponding to
$\lambda = 0$), equation \eqref{iso4_40} gives
\begin{equation}\label{iso100}
\frac{(\rad')^2}{1 + (\rad')^2}\ =\ \alpha^2\ s\ ,
\end{equation}
with
\begin{equation}\label{alpha0}
\alpha\ =\ \frac{2 \lambda}{Q-2}\ .
\end{equation}

 This in turn implies
\begin{equation}\label{E:EL30}
\rad'(s)\ =\ \pm\ \sqrt{\frac{s}{\beta^2 - s}},\qquad \text{where
}\quad \beta\ =\ \frac{1}{\alpha}\ .
\end{equation}

At this point, an observation must be made.  We cannot choose the
sign in \eqref{E:EL30} arbitrarily. In fact, equation
\eqref{iso4_40} implies that $\rad'$ does not change sign on the
interval $[0,R^2/4]$, and one has $\rad'>0$, or $\rad'<0$,
according to whether $\alpha
> 0$ or $\alpha < 0$. On the other hand, if the '$+$' branch of
the square root were chosen in \eqref{E:EL30}, then $\rad$ would
be increasing and, since $\rad \geq 0$ on $(0,R^2/4)$, it would be
thus impossible to fulfill the boundary condition $\rad(R^2/4) =
0$.

These considerations show that it must be $\rad' < 0$ on
$(0,R^2/4)$. We then have to take $\alpha < 0$ (hence $\beta < 0$
as well), and therefore $\lambda < 0$. Equation \eqref{E:EL30}
thus becomes
\begin{equation}\label{E:u'0}
\rad'(s)\ =\ -\ \sqrt{\frac{s}{\beta^2 - s}}\
,\quad\quad\quad\quad 0 \leq s < \frac{R^2}{4}\ .
\end{equation}

We stress that, thanks to the assumption \eqref{l}, and to
\eqref{alpha0}, we have that if \[ 0\ \leq\ s\ < \ \frac{R^2}{4}\
=\ \frac{(Q-2)^2}{4 \lambda^2}\ =\ \frac{1}{\alpha^2}\ =\ \beta^2\
, \]
 then the function $\rad'$ given by \eqref{E:u'0} is smooth on
the interval $[0,R^2/4)$, and satisfies
\[ \underset{s\to \left(\frac{R^2}{4}\right)^-}{\lim}\ \rad'(s)\
=\ -\ \infty\ . \]

 Integrating \eqref{E:u'0} by standard calculus techniques we
find for $s\in [0,R^2/4]$
\begin{align}\label{E:gen_sol0}
\rad(s)\ & =\ \sqrt{s(\beta^2-s)}\ -\
\beta^2\,\tan^{-1}\left(\sqrt{\frac{s}{\beta^2-s}}\right)\ +\ C
\\
& =\ C\ +\  \sqrt{s(\beta^2-s)}\ +\ \beta^2\
\sin^{-1}\left(\frac{\sqrt s}{\beta}\right)\ . \notag
\end{align}

Recalling that $\alpha = \beta^{-1}$, and the equation
\eqref{alpha0}, if we impose the condition $\rad(R^2/4)=0$, we
obtain the solution
\begin{equation}\label{u_bar0}
\rad(s)\ =\  C_{R,\lambda}\ +\ \frac{\sqrt s}{2\lambda}\
\sqrt{(Q-2)^2 - 4 \lambda^2 s}\ +\ \frac{(Q-2)^2}{4\lambda^2}\
\sin^{-1} \left(\frac{2 \ \lambda \sqrt s}{Q-2}\right)\ ,
\end{equation}
where $C_{R,\lambda}$ is given by \eqref{c}. Setting
$u_{R,\lambda}(z) = \rad(|z|^2/4)$, we finally obtain \eqref{uo0}
from \eqref{u_bar0}. We are finally left with proving that such a
$u_{R,\lambda}$ belongs to the class $\mathcal D$. The membership
$u_{R,\lambda}\in \mathcal D$ is equivalent to proving that the
function $s \to \rad(s^2/4)$ is of class $C^1$ in the open interval
$(-R,R)$, and that furthermore $\nabla u_{R,\lambda}\in
C^{0,1}(\Om)$. For the first part, from \eqref{u_bar0} it is clear
that we only need to check the continuity of $\rad'$ at $s=0$. Since
the function is even this amounts to proving that $\rad'(s) \to 0$
as $s\to 0$. But this is obvious in view of \eqref{E:u'0}. Finally,
we have \[ |\nabla u_{R,\lambda}(z) - \nabla u_{R,\lambda}(0)|\ =\
\left|\rad'\left(\frac{|z|^2}{4}\right)\right|\ \leq\ C\ |z|\ , \]
which shows that $\nabla u_{R,\lambda}\in C^{0,1}_{loc}(\Om)$.

\end{proof}

\medskip

In the next Proposition \ref{P:reg1} we complete the analysis of the
regularity of the functions $u_{R,\lambda}$. It suffices to consider
the upper half of the ``normalized'' candidate isoperimetric profile
$E_o\subset \Hn$, where $\p E_o$ is the graph of the function $t =
u_o(z)$, with $u_o = u_{1,\lambda}$ and $\lambda = - (Q-2)$. The
characteristic locus of $E_o$ is given by the two points in $\Hn$
\[
\Sigma\ =\ \{(0,0,\pm \frac{\pi}{8})\}\ .
\]

Unlike its Euclidean counterpart, the hypersurface surface $\p E_o$
is not $C^\infty$ at the characteristic points $(0,0,\pm
\frac{\pi}{8})$.

\medskip

\begin{prop}\label{P:reg1}
The hypersurface $S_o = \p E_o\subset \Hn$ is $C^2$, but not $C^3$,
near its characteristic locus $\Sigma$. However, $S_o$ is $C^\infty$
(in fact, real-analytic) away from $\Sigma$.
\end{prop}

\begin{proof}[\textbf{Proof}]
First, we show that near the characteristic points $(0,0,\pm
\frac{\pi}{8})$ the function $u_o(z)$ given by \eqref{iso} is only
of class $C^2$, but not of class $C^3$. To see this we let
\[
u_1(s) \ =\ \frac{\pi}{8}\ +\ \frac{s}{4} \sqrt{1-s^2}\ -\
\frac{1}{4}\,\sin^{-1}(s)\ , \ 0 \leq s \leq 1\ .
\]
and note that $u_o(z) = u_1(|z|)$ for $0\leq |z| \leq 1$. Therefore,
the regularity of $u_o$ at $|z|=0$ is equivalent to verifying up to
what order of derivatives $n$ one has
\begin{align*}
\lim_{s\to 0^+} u^{(n)}_+(s)\ =\ \lim_{s\to 0^-} u^{(n)}_-(s)
\end{align*}
where $u_+(s) = u_1(s)$ and $u_-(s) = u_1(-s)$.  It is easy to
compute

\begin{align*}
& -u_-'(s) = u_+'(s) = -\frac{1}{2}\,\frac{s^2}{\sqrt{1-s^2}}\ ,
\qquad -u_-''(s) = u_+''(s) = -\frac{1}{2}\frac{s(s^2-2)}{(s^2-1)\sqrt{1-s^2}}\ , \\
& -u_-^{(3)}(s) = u_+^{(3)}(s) =
-\frac{1}{2}\frac{2+s^2}{(s^2-1)^2\sqrt{1-s^2}}\ .
\end{align*}
We clearly have
\begin{align*}
\lim_{s\to 0^-}u_-^{(n)}\ =\ \lim_{s\to 0^+}u_+^{(n)}\quad\text{for
$n =0, 1, 2$}\qquad\text{whereas } \lim_{s\to 0^-}u_-^{(3)}\ =\
1\,\quad\text{and } \lim_{s\to 0^+}u_+^{(3)}\ =\ -1\ .
\end{align*}

This shows the function $t = u_o(z)$ is only $C^2$, but not $C^3$,
near $z=0$.  Next, we investigate the regularity of $\partial E_o$
near $|z| = 1$, that is, at the points where the upper and lower
branches that form $\partial E_o$ meet. To this end, we observe that
$\partial E_o$ can also be generated by rotating around the $t$-axis
the curve in the $(x_1,t)$-plane whose trace is
\[
\{(x_1,t)\,|\,t^2 = u_1(x_1)^2,\quad 0\leq x_1 \leq 1\}\ .
\]

It suffices to show that this curve is smooth ($C^\infty$) across
the $x_1$ axis.  To this end we compute the derivatives of $u_1$. It
is easy to see by induction that for $n \geq 3$
\begin{equation}\label{E:duo}
u_1^{(n)}(x_1)\ =\ (-1)^n\,C_n\,
\frac{P_{n-1}(x_1)}{(x_1^2-1)^{n-1}\sqrt{1 - x_1^2}}\ ,
\end{equation}
where $C_n > 0$ is a constant depending only on $n$, and
$P_{n-1}(x_1)$ is a polynomial in $x_1$ of degree $n-2$. The $n$-th
derivatives of the function $-u_1(x_1)$ clearly takes the same form,
but with a negative sign. Letting $s \to 1^-$ in \eqref{E:duo} we
see that
\[
\frac{d^n}{d x_1^n} u_1, \frac{d^n}{d x_1^n} (-u_1) \longrightarrow
\pm \infty\ , \quad\text{(depending on whether $n$ is odd or even)}\
.
\]

This implies that the curve with equation $t^2 = u_1(x_1)^2$ is
smooth across the $x_1$-axis.

\end{proof}

\medskip

 From Theorem \ref{T:ss0} and Proposition \ref{P:reg1}, we
immediately obtain the following interesting consequence.

\medskip

\begin{thrm}\label{T:ss}
Let $V>0$ be given, and define $R = R(V)>0$ by the formula
\begin{equation}\label{R}
R\ =\ \left(\frac{(Q-2)
\Gamma\left(\frac{Q+2}{2}\right)\,\Gamma\left(\frac{Q-1}{2}\right)}
{\pi^\frac{Q-1}{2}\,\Gamma\left(\frac{Q+1}{2}\right)}\right)^{1/Q}\
V^{1/Q}\ .
\end{equation}
With such choice of $R$, let $\Om = B(0,R) = \{z\in \mathbb R^{2n}
\mid |z|<R\}$. If we take
\begin{equation}\label{lambda}
\lambda\ =\ -\ \frac{Q-2}{R}\ ,
\end{equation}
then the equation \eqref{E:EL-2}, with the Dirichlet condition $u =
0$ on $\p \Om$, admits the spherically symmetric solution $u_R\in
\mathcal D\cap C^2(\Om)$, where
\begin{equation}\label{uo}
u_R(z)\ =\ \frac{\pi R^2}{8}\ +\ \frac{|z|}{4}\ \sqrt{R^2-|z|^2}\
-\ \frac{R^2}{4}\,\sin^{-1}\left(\frac{|z|}{R}\right)\ .
\end{equation}
Furthermore, such $u_R$ satisfies the condition
\begin{equation}\label{vol}
\int_{\Om} u_R(z)\,dz \ =\ \frac{V}{2}\ .
\end{equation}
\end{thrm}

\begin{proof}[\textbf{Proof}]
The first part of the theorem, up to formula \eqref{uo}, is a
direct consequence of Theorem \ref{T:ss0}. We only need to prove
\eqref{vol}. In this respect, keeping in mind the definition
\eqref{R}, it will suffice to prove that
\begin{equation}\label{E:volume}
\int_{\Om} u_R(z)\,dz \ =\
\frac{\pi^\frac{Q-1}{2}\,\Gamma\left(\frac{Q+1}{2}\right)}{2\,(Q-2)
\Gamma\left(\frac{Q+2}{2}\right)\,\Gamma\left(\frac{Q-1}{2}\right)}\
R^Q\ .
\end{equation}

To establish \eqref{E:volume} we note explicitly that $u_R(z) =
\rad(|z|^2/4)$, where
\begin{equation}\label{rad}
\rad(s)\ =\ \frac{\pi R^2}{8}\ +\ \frac{1}{2}\ \sqrt{s(R^2 - 4s)}\
-\ \frac{R^2}{4}\ \sin^{-1} \left(\frac{2 \sqrt s}{R}\right)\ .
\end{equation}

One has therefore
\begin{align*} \int_{\Om} u_R(z)\,dz &\ =\ \int_{|z|<R}
\rad(|z|^2/4)\ dz
\ =\ \sigma_{2n-1}\ \int_0^R \rad(r^2/4)\ r^{2n}\ \frac{dr}{r} \\
&\ =\ 2^{2n-1}\ \sigma_{2n-1}\ \int_0^{\frac{R^2}{4}} \rad(s)\
s^{(Q-4)/2}\ ds\ .
\end{align*}
Integrating by parts the last integral, and using the fact that
$\rad(R^2/4) = 0$, that $\rad$ is smooth at $0$, and \eqref{E:u'0}
(in which now $\beta^2 = \frac{R^2}{4}$), we obtain
\begin{align}\label{E:beta}
\int_{\Om} u_R(z)\,dz \ & =\ \frac{2^{2n}\,\sigma_{2n-1}}{Q-2}\
\int_0^\frac{R^2}{4}
\frac{s^\frac{Q-1}{2}}{\sqrt{\frac{R^2}{4}-s}}\,ds
\\
& =\ \frac{2^{2n+1}\,\sigma_{2n-1}}{Q-2}\ \int_0^\frac{R^2}{4}
s^\frac{Q-2}{2}\,\sqrt{\frac{s}{R^2-4s}}\,ds\ . \notag
\end{align}

With the substitution
\[
t^2\ =\ \frac{R^2 - 4s}{s}, \qquad ds\ =\
\frac{-2R^2t}{(4+t^2)^2}\,\ dt\ ,
\]
the integral \eqref{E:beta} becomes
\begin{align*}
\int_{\Om} u_R(z)\,dz \ & =\ \frac{2^Q\,\sigma_{2n-1}\,R^Q}{Q-2}\
\int_0^\infty \frac{1}{(4+t^2)^\frac{Q+2}{2}}\,dt
\\
& =\ \frac{\sigma_{2n-1}\,R^Q}{4 (Q-2)}\ \int_{\R}
\frac{1}{(1+t^2)^\frac{Q+2}{2}}\,dt\ . \notag \end{align*}

Now we use the formula
\[
\int_{\R} \frac{dt}{(1+t^2)^a}\ =\ \pi^{\frac{1}{2}}\
\frac{\Gamma\left(a - \frac{1}{2}\right)}{\Gamma(a)}\ ,
\]
valid for any $a > 1/2$. We thus obtain
\[
\int_{\Om} u_R(z)\,dz \ =\ \frac{\sigma_{2n-1}\
\pi^{\frac{1}{2}}\,\Gamma\left(\frac{Q+1}{2}\right)}{4\,(Q-2)
\Gamma\left(\frac{Q+2}{2}\right)}\ R^Q\
\]
where $\sigma_{2n-1}$ is the measure of the unit sphere $\mathbb
S^{n-1}$ in $\mathbb R^{2n}$. Finally, using in the latter
equality the fact that
\[
\sigma_{2n-1}\ =\ \frac{2\ \pi^n}{\Gamma(n)}\ =\ \frac{2\
\pi^{\frac{Q-2}{2}}}{\Gamma\left(\frac{Q-2}{2}\right)}\ ,
\]
we obtain \eqref{E:volume}.

\end{proof}

\medskip

With the problem \eqref{E:Gen_Var} in mind, it is convenient to
rephrase part of the conclusion of Theorem \ref{T:ss} in the
following way.

\medskip

\begin{cor}\label{C:R}
Let $V>0$ be given, and for any $R>0$ consider the function $u_R$
defined by \eqref{uo}. There exists $R = R(V)>0$ (the choice of
$R$ is determined by \eqref{R}) such that with $\Om = B(0,R)$ one
has with $u_o = u_R$
\[
\mathcal G[u_o]\ =\ \int_{\Om} u_o(z)\,dz\ =\ \frac{V}{2}\ .
\]

\end{cor}

\medskip

Although the following lemma will not be used until we come to the
proof of Theorem \ref{T:isoine}, it is nonetheless appropriate to
present it at this moment, since it complements Corollary
\ref{C:R}.

\medskip

\begin{lemma}\label{L:ssP}
Let $u_o(z)$ be given by \eqref{uo}, and $\Om = supp(u_o) = B(0,R)$, then
\begin{equation}\label{E:H-perimeter}
\mathcal F[u_o]\ =\ \int_{\Om} \sqrt{|\nabla_z u_o|^2 +
\frac{|z|^2}{4} + <\nabla_z u_o,z^\perp>}\ \ dz\ \ =\
\frac{\pi^\frac{Q-1}{2}\,\Gamma\left(\frac{Q-1}{2}\right)}
{2\Gamma\left(\frac{Q}{2}\right)\,\Gamma\left(\frac{Q-1}{2}\right)}\
R^{Q-1}\ .
\end{equation}
\end{lemma}

\begin{proof}[\textbf{Proof}]
We recall that $u_o(z) = \rad(|z|^2/4)$ where $\rad$ is given by
\eqref{rad}. One has \[
 \nabla_z u_o(z)\ =\
\frac{1}{2}\,\rad'(|z|^2/4)\ z\ ,
\]
and therefore
\[
|\nabla_z u_o(z)|^2\ +\ \frac{|z|^2}{4}\ +\ <\nabla_z
u_o(z),z^\perp> \ =\ \frac{|z|^2}{4}\, \left(1 +
\rad'\left(\frac{|z|^2}{4}\right)^2\right)\ .
\]

We thus obtain
\begin{align*}
& \int_{\Om}  \sqrt{|\nabla_z u_o|^2  + \frac{|z|^2}{4} +
<\nabla_z u_o,z^\perp>}\ dz \ =\ \frac{1}{2}\,
\int_{|z| < R} |z|\,\sqrt{1 + \rad'(|z|^2/4)^2}\,dz \\
&  =\ \frac{\sigma_{2n-1}}{2}\ \int_0^R \sqrt{1 +
\rad'(r^2/4)^{2}}\ r^{2n+1}\ \frac{dr}{r}\ \ =\ 2^{2n-1}\
\sigma_{2n-1}\ \int_0^{R^2/4} \sqrt{1 + \rad'(s)^{2}}\
s^{(Q-3)/2}\ ds\ .
\end{align*}

Formula \eqref{E:u'0}, in which $\beta = -R/2$, gives
\begin{equation*}\label{E:H-per1}
\sqrt{1+\rad'(s)^2}\ =\ \frac{R}{\sqrt{R^2 - 4s}}\ .
\end{equation*}

Inserting this equation in the above integral we obtain
\begin{align*}
\int_{\Om} & \sqrt{|\nabla_z u_o|^2  + \frac{|z|^2}{4} + <\nabla_z
u_o,z^\perp>}\ dz\ \ =\ 2^{2n-1}\ \sigma_{2n-1}\,R\ \int_0^{R^2/4}
s^\frac{Q-4}{2} \sqrt{\frac{s}{R^2 - 4s}}\ ds\ .
\end{align*}

We notice that the last integral above is similar to the one in
\eqref{E:beta}. Proceeding as in the last part of the proof of
Theorem \ref{T:ss}, we finally reach the conclusion.

\end{proof}

\medskip

At this point, recalling that \eqref{E:EL-2} represents the
Euler-Lagrange equation of the unconstrained functional
\eqref{E:unconstrained2}, and keeping \eqref{E:integrants} in
mind, if we combine Theorem \ref{T:ss} with Corollary \ref{C:R},
and take Remark \ref{R:LagrangeM} into account, we obtain the
following result.

\medskip

\begin{thrm}\label{T:crit_pt_obtained}
Let $\mathcal F$ and $\mathcal G$ be the functionals
\[
\mathcal F[u]\ =\ \int_{supp(u)} f(z,\nabla_z u(z))\,dz\ , \qquad
\mathcal G[u]\ =\ \int_{supp(u)} g(u)\, dz\ ,
\]
where $f$ and $g$ are defined in \eqref{E:integrants}. Given $V >
0$, there exists $R = R(V)>0$ (see \eqref{R}) such that the function
$u_o = u_R$ in \eqref{uo} is a critical point in $\mathcal D$ of the
functional $\mathcal F[u]$ subject to the constraint $\mathcal G[u]
= \frac{V}{2}$. This follows from the fact that $u_o$ is a critical
point in $\mathcal D$ of the unconstrained functional $\mathcal
F[u]$ in \eqref{E:unconstrained2}.
\end{thrm}

\medskip

Our next objective is to prove that the function $u_o$ in \eqref{uo}
is: 1) A global minimizer of the variational problem
\eqref{E:Gen_Var};  2) The unique global minimizer. We will need
some basic facts from Calculus of Variations, which we now recall.

\medskip

\begin{dfn}\label{D:Gateaux}
Let $\mathcal V$ be a normed vector space, and $\mathcal D\subset
\mathcal V$. Given a functional $\mathcal F: \mathcal D\to \R$,
$u\in\mathcal D$, and if $\phi$ is $\mathcal D$-admissible at $u$,
one calls
\[
\delta\mathcal F(u;\phi)\ \overset{def}{=}\ \lim_{\epsilon\to
0}\frac{\mathcal F[u + \epsilon \phi] - \mathcal F[u]}{\epsilon}
\]
the G\^ateaux derivative of $\mathcal F$ at $u$ in the direction
$\phi$ if the limit exists.
\end{dfn}

\medskip

\begin{dfn} \label{D:Gateaux-convex}
Let $\mathcal V$ be a normed vector space, and $\mathcal D\subset
\mathcal V$. Consider a functional $\mathcal F: \mathcal D \to
\bar\R$. $\mathcal F$ is said to be convex over $\mathcal D$ if for
every $u\in \mathcal D$, and every $\phi\in \mathcal V$ such that
$\phi$ is $\mathcal D$-admissible at $u$, and $u +\phi \in \mathcal
D$, one has
\[
\mathcal F[u+\phi] - \mathcal F[u]\ \geq\ \delta\mathcal F(u;\phi)\
,
\]
whenever the right-hand side is defined.
We say that $\mathcal F$
is strictly convex if strict inequality holds in the above inequality
except when $\phi\equiv 0$.
\end{dfn}

\medskip

We have the following

\medskip

\begin{thrm}\label{T:unique}
Suppose $\mathcal F$ is convex and proper over a non-empty convex
subset $\mathcal D^*\subset \mathcal V$ (i.e., $\mathcal F
\not\equiv \infty$ over $\mathcal D^*$), and suppose that
$u_o\in\mathcal D^*$ is such that $\delta\mathcal F(u_o;\phi) = 0$
for all $\phi$ which are $\mathcal D^*$-admissible at $u_o$ (that
is, $u_o$ is a critical point of the functional $\mathcal F$), then
$\mathcal F$ has a global minimum in $u_o$. If moreover $\mathcal F$
is strictly convex at $u_o$, then $u_o$ is the unique element in
$\mathcal D^*$ satisfying
\[
\mathcal F[u_o]\ =\ inf\bigl\{\mathcal F[v]\,|\,v\in\mathcal
D^*\bigr\}\ .
\]
\end{thrm}

\begin{proof}[\textbf{Proof}]
Let $u\in \mathcal D^*$, and $u \neq u_o$, then the convexity of
$\mathcal D^*$ implies that $\phi = u - u_o$  is $\mathcal
D^*$-admissible at $u_o$. From Definition \ref{D:Gateaux-convex} we
immediately infer that
\[
\mathcal F[u] - \mathcal F[u_o] \ =\ \mathcal F[u_o + \phi] -
\mathcal F[u_o] \geq\ \delta\mathcal F(u_o;\phi)\ =\ 0\ .
\]
This shows that $\mathcal F$ has a local minimum in $u_o$. When
$\mathcal F$ is strictly convex at $u_o$ we obtain $\mathcal F[u_o +
\phi]
> \mathcal F[u_o]$, for every $\phi\in \mathcal V$ such that $\phi$ is
$\mathcal D^*$-admissible at $u_o$. If $\overline u_o\in \mathcal
D^*$ is another global minimizer of $\mathcal F$, then taking $\phi
= u_o - \overline u_o$, we see that $\mathcal F[u_o] > \mathcal
F[\overline u_o]$. Reversing the roles of $u_o$ and $\overline u_o$
we find $\mathcal F[u_o] = \mathcal F[\overline u_o]$. From the
strict convexity of $\mathcal F$ at $u_o$ we conclude that it must
be $u_o = \overline u_o$.

\end{proof}

\medskip

Our next goal is to adapt the above results to the problem
\eqref{E:Gen_Var}. Given $V>0$ we consider the number $R=R(V)>0$
defined in \eqref{R}, and consider the fixed ball $B(0,R)$. We
consider the normed vector space $\mathcal V(R) = \{u\in C(\overline
B(0,R))\mid  u = 0 \quad \text{on}\quad \p B(0,R)\}$. Let
\begin{align}\label{dstar}
\mathcal D(R)\ =\ & \{u\in \mathcal V(R)\mid u\geq 0\ ,\ u\in
C^2(B(0,R))\cap W^{1,1}(B(0,R))\ ,
\\
& \overline{B}(0,R) = \bigcap \{B(0,R+ \rho)\mid\ supp(u)\subset
B(0,R + \rho)\}\}\ . \notag
\end{align}

We notice that $\mathcal D(R)$ is a non-empty convex subset of
$\mathcal V(R)$, and that for every $u\in \mathcal D(R)$ one has $u
=0$ on $\p B(0,R)$. Let $h=h(z,u, p)$ be the function in
\eqref{E:integrants} and consider the functional
\eqref{E:unconstrained2}. Given $u\in \mathcal D(R)$ and $\phi$
which is $\mathcal D(R)$-admissible at $u$, in view of Theorem
\ref{T:balogh}, we see that $\mathcal F$ is G\^ateaux differentiable
at $u$ in the direction of $\phi$, and
\begin{align}\label{E:Gateaux2}
\delta \mathcal F(u;\phi)\ & =\ \int_{B(0,R)}
\bigg\{h_u(z,u(z),\nabla u(z))\,\phi(z) \ +\ <\nabla_p
h(z,u(z),\nabla u(z)),\nabla \phi(z)>\bigg\}\,dz
\\
& =\ \int_{B(0,R)} \left\{\frac{<\nabla_z u + z^\perp/2,
\nabla_z\phi>}{| \nabla_z u + z^\perp/2|} \ +\
\lambda\,\phi\right\}\,dz\ , \notag
\end{align}
where in the above $<\cdot,\cdot>$ denotes the standard inner
product on $\R^{2n}$. One has the following well-known sufficient
condition for the convexity (strict convexity) of $\mathcal F$.

\medskip

\begin{prop}\label{P:s-convex}
If for $a.e.\ z\in B(0,R)$, for all $u\in\mathcal D(R)$ and $p =
\nabla u$, the function $h$ in the definition of $\mathcal F$
satisfies for every $\phi$ which is $\mathcal D(R)$-admissible at
$u$, and every $q= \nabla \phi$,
\begin{equation}\label{E:S-convex}
h(z,u+v, p+ \phi)\ -\ h(z,u, p) \ \geq\ h_u(z,u, p)\,\phi \ +\
<\nabla_{ p} h(z,u, p), q>\,
\end{equation}
then $\mathcal F$ is convex on $\mathcal D(R)$. If, instead, the
strict inequality holds unless $v=0$ and $ q= 0$, then $\mathcal F$
is strictly convex.
\end{prop}

\begin{proof}[\textbf{Proof}]
Let $u\in \mathcal D(R)$, and let $\phi$ be $\mathcal
D(R)$-admissible at $u$. Using \eqref{E:Gateaux2}, we obtain
\begin{align*} \mathcal F[u + \phi] -
\mathcal F[u] &\ =\ \int_{B(0,R)} \big\{h(z,u(z)+\phi(z),\nabla u(z)
+ \nabla \phi(z)) -
h(z,u(z),\nabla u(z))\big\}\,dz \\
 & \ \geq \ \int_{B(0,R)} \big\{h_u(z,u(z),\nabla u(z))\,\phi(z) \ +\
<\nabla_{p}h(z,u(z),\nabla u(z)),\nabla\phi(z)>\big\}\,dz\\
 &\ =\ \delta\mathcal F(u;\phi)\ .
\end{align*}

Appealing to Definition \ref{D:Gateaux-convex} the conclusion
follows.

\end{proof}

\bigskip

Our next goal is to prove that the unconstrained functional
$\mathcal F$ in \eqref{E:unconstrained2} is convex on the convex set
$\mathcal D(R)$. Since each one of them has an independent interest,
we will provide two different proofs of this fact. The former is
based on the following linear algebra lemma, which is probably
well-known, and whose proof we have provided for the reader's
convenience.

\medskip

\begin{lemma}\label{L:pos_semi_def}
Let $\mathbf A=[A_{ij}]$ be an $m \times m$ matrix with entries
given by
\[
A_{ij} \ =\ \delta_{ij} - \frac{a_i\,a_j}{D} \qquad \text{where } D\
=\ \sum_{i=1}^m a_i^2\ \not=\ 0\ ,
\]
then $\mathbf A$ has $\lambda = 0$ as an eigenvalue of multiplicity
one, and $\lambda = 1$ as an eigenvalue of multiplicity $m-1$.
\end{lemma}

\begin{proof}[\textbf{Proof}]
First, consider the matrix $\mathbf I - \mathbf A$, which takes the
form
\[
\frac{1}{D}\
\begin{pmatrix}
a_1\,a_1 & a_1\,a_2 & a_1\,a_3 & \cdots & a_1\,a_m \\
a_2\,a_1 & a_2\,a_2 & a_2\,a_3 & \cdots & a_2\,a_m \\
\vdots   & \cdots & \vdots  &  \ddots & \vdots \\
a_m\,a_1 & a_m\,a_2 & a_m\,a_3 & \cdots & a_m\,a_m
\end{pmatrix}\ .
\]
It is easy to see that an equivalent row-echelon form of the matrix
has the last $m-1$ rows containing all zeros, thus $\mathbf I -
\mathbf A$ is a matrix of rank one. From the rank-nullity theorem we
conclude that $\lambda =1$ is an eigenvalue of $\mathbf A$ of
multiplicity $m-1$. We are thus left with showing the $\lambda = 0$
is a simple eigenvalue. For this we show that $det(\mathbf A) = 0$.
Observe that $det(\mathbf A) = D^{-m} det(\mathbf B)$, where
\[
\mathbf B \ =\
\begin{pmatrix}
D-a_1^2 & -a_1\,a_2 & -a_1\,a_3 & \cdots & -a_1\,a_m \\
-a_1\,a_2 & D - a_2^2 & -a_2\,a_3 & \cdots & -a_2\,a_m \\
\vdots &    \cdots    & \vdots    & \ddots & \vdots \\
-a_m\,a_1 & -a_m\,a_2 & -a_m\,a_3 & \cdots & D - a_m^2
\end{pmatrix}\ .
\]

To continue the computation of $det(\mathbf B)$, we replace rows
$R_j$ by $a_1\,R_j - a_j \,R_1$ for $j = 2,...,m$ and observe that
$a_1\,R_j - a_j \,R_1$ takes the form
\[
a_1\,R_j \ -\ a_j\,R_1 = [-a_j\,D \ 0\  \cdots\ 0\ a_1\,D\ 0\
\cdots\ 0]\ .
\]

We then have
\[
det(\mathbf B) \ =\det(\mathbf C)\ ,
\]
where
\[
\mathbf C\ =\
\begin{pmatrix}
D - a_1^2  &  -a_1\,a_2 & \cdots  &   \cdots  & \cdots & -a_1\,a_m \\
-a_2\,D    &  a_1\,D    &   0     &   \cdots  & \cdots &     0     \\
-a_3\,D    &    0       & a_1\,D  &      0    & \cdots &     0     \\
\vdots     &  \cdots    & \cdots  &   \cdots  & \ddots &   \vdots  \\
-a_m\,D    &    0       &   0     &      0    & \cdots &    a_1\,D
\end{pmatrix}\ .
\]

To compute $det(\mathbf C)$ we take advantage of the special
structure of the matrix, and consider
\[
\mathbf C\,\mathbf C^T\ =\ D^2
\begin{pmatrix}
-a_1 & -a_2 & -a_3 & \cdots & -a_m\\
-a_2 & a_2^2 & a_2\,a_3 & \cdots & a_2\,a_m\\
\vdots & \cdots & \cdots & \ddots & \vdots \\
-a_m & a_m\,a_2 & a_m\,a_3 & \cdots & a_m^2\\
\end{pmatrix}\ .
\]

We note that if either $a_2 = 0$ or $a_3 = 0$, then the matrix
$\mathbf C \mathbf C^T$ has a column of zeros, and therefore its
determinant vanishes. Suppose then that $a_2, a_3 \not= 0$.
Replacing rows $R_2$ and $R_3$ by $R_2 + a_2\,R_1$ and $R_3 +
a_3\,R_1$ respectively, we see that the new rows two and three have
first entries given by $-a_2 - a_1\,a_2$ and $-a_3 - a_1\,a_2$,
whereas all the remaining entries vanish. Either one of these rows
is already a zero row or else, using one to eliminate the other, we
obtain a row of zeros, and therefore we conclude that $det(\mathbf
C\mathbf C^T) = 0$. Hence, $det(\mathbf A) = D^{-m}\,det(\mathbf B)
= D^{-m}\,det(\mathbf C) = 0$.  This completes the proof of the
lemma.

\end{proof}

\medskip

\medskip

\begin{prop}\label{P:convex}
Given $V>0$, let $R = R(V)>0$ be as in \eqref{R}. The functional
$\mathcal F$ in \eqref{E:unconstrained2} is convex on $\mathcal
D(R)$. As a consequence, the function $u_R$ in \eqref{uo} is a
global minimizer of $\mathcal F$ on $\mathcal D(R)$.
\end{prop}

\begin{proof}[\textbf{Proof}]
Considering the integrand $h(z,u,p) = |p + \frac{z^\perp}{2}| +
\lambda u$ in the functional $\mathcal F$ in
\eqref{E:unconstrained2}, we have
\[
\begin{cases}
h_{p_i} = (p_i + \frac{z^\perp_i}{2})/|p + \frac{z^\perp}{2}|\ , \\
h_{p_i\,p_j} \ =\ \frac{1}{|p + \frac{z^\perp}{2}|}
\left\{\delta_{ij}
- \frac{(p_i + \frac{z^\perp_i}{2})\,(p_j + \frac{z^\perp_j}{2})}{|p + \frac{z^\perp}{2}|^2} \right\}\ ,\\
h_{u,p_i}\ = \ 0\ ,
\end{cases}
\]
where in the above we have let
\begin{equation}\label{perp}
z^\perp_i \ =\
\begin{cases}
\ \,\, y_i \quad\text{if } 1 \leq i \leq n \\
-x_i \quad\text{if } n+1 \leq i \leq 2n
\end{cases}
\end{equation}
The hessian of $h$ with respect to the variable $(u,p)\in \R \times
\R^{2n}$ now takes the form
\[
\nabla^2 h(u,p) \ =\
\begin{pmatrix}
0 & 0 & \cdots & 0 \\
0 &   &   &   &   \\
\vdots & & \mathcal A  & \\
0   &   &   &   \\
\end{pmatrix}
\]
where, aside from the multiplicative factor $1/|p + z^\perp/2|$, the
block $\mathcal A$ takes the form of the matrix $\mathbf A$ in Lemma
\ref{L:pos_semi_def}. We thus conclude that the eigenvalues of
$\nabla^2 h(u,p)$ are $\lambda = 0$ (of multiplicity two) and
$\lambda = 1/|p + z^\perp/2|$ of multiplicity $2n-1$. Thus, from
Theorem \ref{T:balogh}, for a.e. $z\in  B(0,R)$, the function $(u,p)
\to h(z,u,p)$ is convex.  This in turn implies that $\mathcal F$ is
convex on $\mathcal D(R)$. From Theorems \ref{T:crit_pt_obtained}
and \ref{T:unique} we conclude that $u_R$ is a global minimizer of
$\mathcal F$ on $\mathcal D(R)$.

\end{proof}

\medskip

We next prove a slightly stronger result than Proposition
\ref{P:convex}, namely the convexity of the function in $\R^{2n}$
which defines the integrand in $\mathcal F$ in
\eqref{E:unconstrained2}. The proof of this result is based on the
following lemma.

\medskip

\begin{lemma}\label{L:convexity}
Let $\alpha \in \R^{2n}$ be fixed, with $\alpha \not= 0$, then one
has
\[
f(q)\ \overset{def}{=}\ |\alpha|\ |q|^2\ -\ (|q+\alpha| - |\alpha|)
<q,\alpha>\ \geq \ 0\ ,\quad\quad\quad \text{for every}\quad q\in
\R^{2n}\ .
\]
\end{lemma}

\begin{proof}[\textbf{Proof}]
We observe that $f(0) = f(-\alpha) = 0$, and that $f\in
C^\infty(\R^{2n}\setminus \{-\alpha\})$. We want to analyze the
possible critical points in $\R^{2n}\setminus \{-\alpha\}$ of the
function $f$. It is easier to reduce the problem by introducing
spherical coordinates. Let $\alpha = r_0 \omega_0$, with
$\omega_0\in \mathbb S^{2n-1}$ and $r_0>0$, we can consider a system
of spherical coordinates in which the ``north pole" coincides with
$\omega_0$ and the colatitude angle $\theta$ denotes the angle
formed by the vector $q\in \R^{2n}\setminus \{-\alpha\}$ with
$\omega_0$. In such a system we let $q= r \omega$, with $\omega\in
\mathbb S^{2n-1}$, and $r = |q|$, so that $cos\ \theta =
<\omega,\omega_0>$. We observe that the function $f$ is constant on
every $2n-2$ dimensional sub-sphere $\sin \theta = const$ of the
unit sphere $\mathbb S^{2n-1}\subset \R^{2n}$, and we want to
exploit these symmetries of $f$. For $z\in \R^{2n}$ we let $r = r(z)
= |z|$, and $\theta = \theta(z) = \cos^{-1}(<z/r,\omega_0>)$.
Writing $f(z) = f(r(z),\theta(z))$, we are thus led to consider
\[
f(r,\theta)\ =\ f(q)\ =\ r_0\ r^2  - \big(\sqrt{r_0^2 + r^2 + 2 r_0
r \cos \theta}\ -\ r_0\big) r_0\ r\ \cos\ \theta\ .
\]

If we now set $t = r/r_0$, then we can consider the function
\[
g(t,\theta)\ =\ \frac{1}{r_0^3}\ f(r_0 t,\theta)\ =\ t^2 -
\big(\sqrt{1 + t^2 + 2 t \cos \theta}\ -\ 1\big)\ t \ \cos \theta\ ,
\]
for $(t,\theta)\in Q =  [0,\infty) \times [0,\pi]$, with $(t,\theta)
\not= (1,\pi)$. When $t=0$, then $q= 0$ and we have already observed
that $f(0) = 0$. When $\theta = 0$, then $q = \rho \alpha$ for some
$\rho \geq 0$, one readily recognizes that $f(\rho \alpha) = 0$.
Finally, when $t\geq 0$ and $\theta = \pi$ we have $g(t,\pi) = 0$ if
$0\leq t \leq 1$, and $g(t,\pi)>0$ for $t>1$. In conclusion, we have
$f(q) = 0$ for $q = \rho \alpha$ for some $\rho \geq - 1$, whereas
we have $f(q)>0$ for $q = \rho \alpha$ with $\rho < - 1$. We now
consider the possible critical points of $f$. Using the chain rule
we see that
\[
\nabla f\ =\ \frac{f_r}{r}\ z\ -\ \frac{f_\theta}{r \sin \theta}\
\left(\omega_0 - \frac{\cos \theta}{r}\ z\right)\ . \]

Since $<z,\omega_0 - \frac{\cos \theta}{r} z> = 0$, we find \[
|\nabla f|^2\ =\ f_r^2 \ +\ \frac{f_\theta^2}{r^2 \sin^2 \theta}
\left|\omega_0 - \frac{\cos \theta}{r} z\right|^2\ =\ f_r^2\ +\
\frac{1}{r^2}\ f^2_\theta\ ,
\]
which allows us to conclude that $\nabla f$ vanishes outside of the
set of points $q = \rho \alpha$ with $\rho \geq - 1$, if and only if
$f_r= f_\theta = 0$ at interior points of $Q$. This is equivalent to
studying the interior critical points of the function $g(t,\theta)$
in $Q$. One has
\begin{align}
\nabla g(t,\theta)\  =\ & \bigg(2t - \big(\sqrt{1 + t^2 + 2 t \cos
\theta}\ -\ 1\big)\ \cos \theta - \frac{(t + \cos \theta) t \cos
\theta}{\sqrt{1 + t^2 + 2 t \cos \theta}}\ ,
\\
& \big(\sqrt{1 + t^2 + 2 t \cos \theta}\ -\ 1\big)\ t \sin \theta +
\frac{t^2 \sin \theta\ \cos \theta}{\sqrt{1 + t^2 + 2 t \cos
\theta}}\bigg) \notag\\
& =\ (g_t,g_\theta)\ . \notag
\end{align}

Since now $0<\theta < \pi$ it is clear that $g_\theta = 0$ if and
only if
\begin{equation}\label{zerogteta}
\big(\sqrt{1 + t^2 + 2 t \cos\ \theta}\ -\ 1\big)\ =\ -\ \frac{t
\cos\ \theta}{\big(\sqrt{1 + t^2 + 2 t \cos\ \theta}}\ .
\end{equation}

On the other hand, we see that $g_t = 0$ at points where
\eqref{zerogteta} holds if and only if
\[
g_t\ =\ 2t -\ \frac{t^2 \cos \theta}{\sqrt{1 + t^2 + 2 t \cos\
\theta}}\ =\ 0\ ,
\]
which is equivalent to
\begin{equation}\label{crit}
2\ =\ \frac{t \cos \theta}{\sqrt{1 + t^2 + 2 t \cos\ \theta}}\ .
\end{equation}

It is clear that if $\pi/2 < \theta < \pi$, then \eqref{crit} has no
solutions. Suppose then that $0<\theta < \pi/2$. In this range,
equation \eqref{crit} is equivalent to
\[
4\ =\ \frac{t^2 \cos^2 \theta}{1 + t^2 + 2 t \cos\ \theta}\ ,
\]
which is in turn equivalent to \[ (4 - \cos^2 \theta) t^2 + 8 t \cos
\theta + 4\ =\ 0\ .
\]

An easy verification which we leave to the reader shows that the
latter equation has no solutions $t>0$ in the range $0<\theta <
\pi/2$. In conclusion, the function $g(t,\theta)$, and therefore has
no interior critical points. Therefore, $g(t,\theta) \geq 0$ for
every $(t,\theta)\in Q$. This allows to conclude that $f(q) \geq 0$
for all $q\in \R^{2n}$, thus completing the proof of the lemma.

\end{proof}

\medskip

At this point we observe that Lemma \ref{L:convexity} provides an
alternative proof of Proposition \ref{P:convex}. It suffices in fact
to consider for every $u\in \mathcal D(R)$ and every $\phi$ which is
$\mathcal D(R)$-admissible at $u$, the vectors $\alpha(z) = \nabla
u(z) + z^\perp/2$, $q(z) =\nabla \phi(z)$. Let $\mathcal F$ be given
by \eqref{E:unconstrained2} and recall \eqref{E:Gateaux2}. One has,
\begin{align}\label{E:diff}
\mathcal F[u + \phi]\ & -\ \mathcal F[u] \\
\notag &\ =\ \int_{B(0,R)} \big\{|\nabla_z u + z^\perp/2 + \nabla_z
\phi | \ -\ |\nabla_z u + z^\perp/2 | \ +\
\lambda\, \phi\big\} \,dz \\
\notag &\ =\ \int_{B(0,R)} \left\{\frac{2\,<\nabla_z\phi,\nabla_z u
+ z^\perp/2>\ +\ |\nabla_z \phi |^2}{|\nabla_z u + z^\perp/2 |\ +\
|\nabla_z u + z^\perp/2 + \nabla_z\phi| } \ +\ \lambda\, \phi
\right\}\,dz\ . \notag
\end{align}

From Theorem \ref{T:balogh} we know that there exists $Z\subset
\Om$, with $|\Om \setminus Z| = 0$, such that $|\alpha(z)| \not= 0$
for every $z\in Z$. We intend to show that for every $z\in Z$ we
have
\begin{equation}\label{E:ge}
\frac{2< q,\alpha> + |q|^2}{|\alpha| + |\alpha + q|}\ \geq\
\frac{<q,\alpha>}{|\alpha|}\ .
\end{equation}

This would imply
\begin{equation}\label{si}
\frac{2<\nabla_z\phi,\nabla_z u + z^\perp/2>\ +\ |\nabla_z \phi
|^2}{|\nabla_z u + z^\perp/2 |\ +\ |\nabla_z u + z^\perp/2 +
\nabla_z\phi| }\  \geq\ \frac{<\nabla_z\phi,\nabla_z u +
z^\perp/2>}{| \nabla_z u + z^\perp/2|}\ ,
\end{equation}
which would prove that $\mathcal F$ is convex. For every $z\in Z$
the inequality \eqref{E:ge} is easily seen to be equivalent to
\begin{equation}\label{q}
(|q + \alpha|\ -\ |\alpha|)\ <q,\alpha> \ \leq \ | q|^2\,|\alpha| \
,
\end{equation}
which is true in view of Lemma \ref{L:convexity}. Finally, we give
the proof of Theorem \ref{T:isoprofile}.

\medskip

\begin{proof}[\textbf{Proof of Theorem \ref{T:isoprofile}}]
We fix $V>0$ and consider the collection of all sets $E\in\mathcal
E$ such that $V=|E|$. We want to show that the problem of minimizing
$P_H(E;\Hn)$ within this subclass admits a unique solution, and that
the latter is given by \eqref{uo}, in which the parameter $R = R(V)$
has been chosen as in \eqref{R}. According to condition (i) in the
definition of the class $\mathcal E$, we have $V/2 = |E\cap \Hn_+|$.
Still from assumption $(i)$, and in view of Theorem
\ref{T:iso-symm}, it is enough to minimize
$P_H(E;\overline{\Hn_+})$. This is an important point. In fact,
Theorem  \ref{T:iso-symm} states that, if $E$ is an isoperimetric
set, i.e., if $E$ minimizes $P_H(\circ;\Hn)$ under the constraint
$|E|= V$, then \begin{equation}\label{p} P_H(E;\overline{\Hn_+})\ =\
P_H(E;\overline{\Hn_-})\ .
\end{equation}

This implies that the minimizer must be sought for within the class
of sets $E\in \mathcal E$ such that $|E| = V$, and for which
\eqref{p} holds, which is in turn equivalent to proving existence
and uniqueness of a global minimizer in the class $\mathcal D(R)$
defined by \eqref{dstar}. The existence of a global minimizer
follows from Proposition \ref{P:convex}, and such global minimizer
is provided by the spherically symmetric function $u_R$ in
\eqref{uo}. We are thus left with proving its uniqueness. The latter
will follow if we can prove that for every $\mathcal
D(R)$-admissible function $\phi$ at $u_R$ the strict inequality
\[
\mathcal F[u_R + \phi]\ >\ \mathcal F[u_R]
\]
holds, unless $\phi \equiv 0$. This will follow from the strict
inequality in \eqref{si} for every $z\in Z$, with $u$ replaced by
the function $u_R$ in \eqref{uo}, unless $\phi \equiv 0$ in
$B(0,R)$. Such strict inequality is equivalent to proving strict
inequality in \eqref{q} on the set $Z$, with $q(z) = \nabla \phi(z)$
and $\alpha(z) = \nabla u_R(z) + z^\perp/2$. We emphasize here that,
in view of \eqref{uo}, the vector-valued function $\alpha(z)$ only
vanishes at $z = 0$. Keeping in mind that $u_R\in C^2(B(0,R))$, and
that, since $\phi$ is $\mathcal D(R)$-admissible at $u_R$, we have
$\phi\in C^2(B(0,R))$, and $\phi = 0$ on $\p B(0,R)$, an analysis of
the proof of Lemma \ref{L:convexity}, brings to the conclusion that
the desired strict inequality holds, unless either $\nabla
\phi\equiv0$, in which case we conclude $\phi \equiv 0$, or there
exists a function $\rho\in C^1(B(0,R))$, with $\rho \geq - 1$,  and
such that for every $z\in Z$
\begin{equation}\label{grad}
\nabla \phi(z)\ =\ \rho(z)\ \left(\nabla u_R(z) +
\frac{z^\perp}{2}\right)\ .
\end{equation}

We remark explicitly that the possibility $\rho \equiv const$ in
\eqref{grad} is forbidden by the fact that the vector field $z\to
\nabla u_R(z) + z^\perp/2$ is not conservative in $B(0,R)$.
Furthermore, since the functions in both sides of \eqref{grad} are
in $C^1(B(0,R))$, the validity of the inequality for every $z\in Z$
is equivalent to its being valid on the whole $B(0,R)$.

We thus want to show that \eqref{grad} cannot occur. To illustrate
the idea, we focus on the case $n=1$ and leave the trivial
modifications to the interested reader. We argue by contradiction
and suppose that \eqref{grad} hold. This means
\[
\phi_x\ =\ \rho \left(u_{R,x} + \frac{y}{2}\right)\ , \ \phi_y\ =\
\rho \left(u_{R,y} - \frac{x}{2}\right)\ .
\]

Since $\phi \in C^2(B(0,R))$, differentiating the first equation
with respect to $y$ and the second with respect to $x$, and keeping
in mind that $u_R$ is spherically symmetric (see \eqref{uo}), from
the fact that $\phi\in C^2(B(0,R))$, and therefore $\phi_{xy}
=\phi_{yx}$, we infer that we must have
\begin{equation}\label{tr}
\left(\frac{x}{2} - \overline u'\ \frac{y}{2}\right)\ \rho_x\ +\
\left(\frac{y}{2} + \overline u'\ \frac{x}{2}\right)\ \rho_y \ +\
\rho\ =\ 0\ ,
\end{equation}
where, we recall, $u_R(z) = \rad(|z|^2/4)$, see \eqref{rad}. We now
fix a point $z_0\in B(0,R)\setminus\{0\}$, and consider the
characteristic curve starting at $z_0 = (x_0,y_0)$, $z(s) =
z(s,z_0)$ of the transport equation \eqref{tr}. Letting $z(s) =
(x(s),y(s))$, we know that such curve satisfies the system
\begin{equation}\label{system}
\begin{cases}
x'\ =\ \frac{x}{2} - \overline u'\ \frac{y}{2}\ ,\quad\quad x(0)\ =\
x_0\ ,
\\
y'\ =\ \frac{y}{2} + \overline u'\ \frac{x}{2}\ ,\quad\quad y(0)\ =\
y_0\ .
\end{cases}
\end{equation}

It is clear that $s\to \rho(z(s))$ satisfies the Cauchy problem
\[
\frac{d}{ds}\ \rho(z(s))\ =\ -\ \rho(z(s))\ ,\quad\quad\quad
\rho(z(0)) = \rho(z_0)\ ,
\]
and therefore
\begin{equation}\label{ro}
\rho(z(s))\ =\ \rho(z(s,z_0))\ =\ \rho(z_0)\ e^{-s}\ .
\end{equation}

Multiplying the first equation in \eqref{system} by $x$, and the
second by $y$, we find
\[
\frac{d}{ds} |z(s)|^2\ =\ |z(s)|^2\ ,
\]
which gives
\begin{equation}\label{z}
|z(s)|^2\ =\ |z_0|^2\ e^s\ .
\end{equation}

It is clear that $-\infty < s \leq 2 \log(R/|z_0|)$. For every $s$
in this range, we obtain from \eqref{grad}, \eqref{ro},  and from
\eqref{E:u'0},
\[
\nabla \phi(z(s))\ =\ \frac{\rho(z_0) e^{-s}}{2} \left(-
\frac{|z(s)|}{\sqrt{R^2 - |z(s)|^2}}\ z(s)\ +\ z(s)^\perp\right)\ .
\]

Using \eqref{z}, we finally obtain \[ |\nabla \phi(z(s))|^2\ =\
\frac{\rho(z_0)^2 e^{-2s}}{4} |z_0|^2 e^{2s}\
\left[\frac{|z(s)|^2}{R^2 -|z(s)|^2}\ +\ 1\right]\  .
\]

Letting $s\to - \infty$ in the latter equation, we reach the
conclusion
\[
|\nabla \phi(0)|^2\ =\ \frac{\rho(z_0)^2 |z_0|^2}{4}\ ,
\]
which contradicts the continuity of $|\nabla\phi|$ at $z=0$, unless
$\rho\equiv 0$. But this would contradict our assumptions on $\rho$.
We conclude that $u_R$ given by \eqref{uo} is the unique minimizer
to the variational problem \eqref{E:Gen_Var} in $\mathcal D(R)$.

\end{proof}

\medskip

\begin{rmrk}\label{R:chmy}
We mention that an alternative proof of the uniqueness of the global
minimizer $u_R$ in Theorem \ref{T:isoprofile} could be obtained by
the interesting comparison Theorem C' on p.163 in \cite{CHMY}.
\end{rmrk}

\medskip

\begin{prop}\label{P:necessary}
Suppose $E\in\mathcal E$ is a critical point of the $H$-perimeter
subject to the constraint $|E| = const$, then $S =
\partial E$ has constant $H$-mean curvature. In particular, the
isoperimetric set $E_o$ found in Theorem \ref{T:isoprofile} is a set
of constant positive $H$-mean curvature $\mathcal H =
\frac{Q-2}{R}$.
\end{prop}

\begin{proof}[\textbf{Proof}]
Let $E\in\mathcal E$ be given and let $u$ be the function
describing $\partial E$ in $\Hn_+$. To prove that $\partial E$ has
constant $H$-mean curvature we could  appeal to Remark
\ref{R:interpretation}. Instead, we proceed directly as follows.
We recall that $u(z) = \rad(|z|^2/4)$ for some $C^2$ function
$\rad$, and the assumptions that $E$ is a critical point of the
$H$-perimeter means that $\rad$ satisfies \eqref{iso3_10}. From
the discussion in the proof of Theorem \ref{T:ss0}, the left hand
side of \eqref{iso3_10} (that is the Euler-Lagrange equation)
becomes
\[
r\,F'(r)\ +\ (Q-2) F(r)
\]
where $F(r)$ is given by \eqref{F0}.  A simple computation gives
\[
F'(r)\ =\ \frac{r^2\,\rad''(r^2/4)\ -\ 2\,\rad'(r^2/4)\,(1 +
\rad'(r^2/4)^2)}{2\,r^2\,(1 + \rad'(r^2/4)^2)^\frac{3}{2}}\ ,
\]
and therefore we have
\begin{equation}\label{EL=Hcurv}
r\,F'(r)\ +\ (Q-2)\,F(r)
\ =\
\frac{2\,(Q-3)\rad'(r^2/4)(1 + \rad'(r^2/4)^2)\ +\ r^2 \rad''(r^2/4)}
{2\,r\,(1 + \rad'(r^2/4)^2)^\frac{3}{2}}\ .
\end{equation}

Rewriting the Euler-Lagrange equation \eqref{iso3_10} for such functions $u$ (or $\rad$) we have
\begin{equation}\label{EL-cylindrical}
\frac{2\,(Q-3)\rad'(r^2/4)(1 + \rad'(r^2/4)^2)\ +\ r^2 \rad''(r^2/4)}
{2\,r\,(1 + \rad'(r^2/4)^2)^\frac{3}{2}}\ =\ \lambda
\end{equation}
where $\lambda$ is of course a constant.  We make a change of notation by
letting $s = r^2/4$ in \eqref{EL-cylindrical}, we found

\begin{equation}\label{EL-cylindrical2}
\frac{(Q-3)\rad'(s)(1 + \rad'(s)^2)\ +\ 2\,s\, \rad''(s)}
{2\,\sqrt{s}\,(1 + \rad'(s)^2)^\frac{3}{2}}\ =\ \lambda\ .
\end{equation}
Comparing \eqref{EL-cylindrical2} with \eqref{tphi}, we infer that
the $H$-mean curvature of such surfaces is
\[
\mathcal H\ =\ -\, \frac{(Q-3)\rad'(s)(1 + \rad'(s)^2)\ +\ 2 s\,
\rad''(s)} {2\,\sqrt{s}\,(1 + \rad'(s)^2)^\frac{3}{2}}\ =\
-\,\lambda\ .
\]
If the set $E_o$ is described by $u_R(z)$, where $u_R(z)$ is given
by \eqref{uo}, then from \eqref{lambda} in Theorem \ref{T:ss} we
conclude that the $H$-mean curvature of $E_o$ is given by
\[
\mathcal H \ =\ \frac{Q-2}{R}\ .
\]

\end{proof}

This completes proof of Theorem \ref{T:isoprofile}.

\medskip

\begin{proof}[\textbf{Proof of Theorem \ref{T:isoine}}]
We have already established the restricted isoperimetric
inequality. Furthermore, the invariance of the isoperimetric
quotient with respect to the group translations and dilations is a
consequence of Propositions \ref{P:invariance} and
\ref{P:invariance2}. We are left with the computation of the
constant $C_Q$. To this end, we use the set $E_R$ described by
$u_o$. We note that the integrals \eqref{E:volume} and
\eqref{E:H-perimeter} give $|E_R|/2$ and $P_H(E;\Hn_+)$
respectively, and therefore after some elementary simplifications
we obtain
\[
C_Q\ =\ \frac{|E_R|^\frac{Q-1}{Q}}{P_H(E_R;\Hn)} \ =\
\frac{(Q-1)\Gamma\left(\frac{Q}{2}\right)^\frac{2}{Q}}{Q^\frac{Q-1}{Q}(Q-2)
\Gamma\left(\frac{Q+1}{2}\right)^\frac{1}{Q}\pi^\frac{Q-1}{2Q}}\ .
\]

This completes the proof.

\end{proof}


\vskip 0.6in


\begin{thebibliography}{99}

\bibitem[A]{A}
A. D. Alexandrov, \emph{A characteristic property of spheres},
Ann. Mat. Pura Appl., (4) \textbf{58}~1962, 303-315.



\bibitem[B]{B}
Z. M. Balogh, \emph{Size of characteristic sets and functions with
prescribed gradients}, J. Reine Angew. Math., \textbf{564}~2003,
63-83.

\bibitem[BSV]{BSV}
V. Barone Adesi, F. Serra Cassano and D. Vittone, \emph{The
Bernstein problem for intrinsic graphs in the Heisenberg group and
calibrations}, preprint, 2006.

\bibitem[Be]{Be}
  A. Bella\"{\i}che, \emph{The tangent space in sub-Riemannian geometry.
Sub-Riemannian geometry,}, Progr. Math., \textbf{144}~(1996), Birkh\"auser, 1-78.


\bibitem[BM]{BM}
M. Biroli \& U. Mosco, \emph{Sobolev and isoperimetric inequalities for Dirichlet forms on homogeneous spaces}, Pot. Anal., \textbf{4}~(1995), 311-324.

\bibitem[BC]{BC}
M. Bonk \& L. Capogna, \emph{Mean Curvature flow  and the
isoperimetric profile of the Heisenberg group}, preprint, 2005.



\bibitem[CDG]{CDG}
L. Capogna, D. Danielli \& N. Garofalo, \emph{The geometric Sobolev embedding for vector fields and the isoperimetric inequality}, Comm. Anal. and Geom., \textbf{2}~(1994), 201-215.



\bibitem[CG]{CG}
L. Capogna \& N. Garofalo, \emph{Ahlfors type estimates for
perimeter measures in Carnot-Carath\'eodory spaces}, J. Geom.
Anal., to appear.


\bibitem[Ca]{Ca}
C. Carath\'eodory, \emph{Untersuchungen \"uber die Grundlangen der Thermodynamik}, Math. Ann., \textbf{67}~(1909), 355-386.

\bibitem[CH]{CH}
J.H.Cheng \& J.F. Hwang, \emph{Properly embedded and immersed
minimal surfaces in the Heisenberg group}, Bull. Austral. Math.
Soc., \textbf{70}~(2004), no. 3, 507-520.


\bibitem[CHMY]{CHMY}
J.H. Cheng, J. F. Hwang, A. Malchiodi \& P. Yang, \emph{Minimal
surfaces in pseudohermitian geometry and the Bernstein problem in
the Heisenberg group}, revised version 2004, Ann. Sc. Norm. Sup.
Pisa, \textbf{1}~(2005), 129-177.

\bibitem[Ch]{Ch}
W. L. Chow, \emph{\"Uber Systeme von linearen partiellen Differentialgleichungen erster Ordnung}, Math. Annalen, \textbf{117} (1939), 98-105.




\bibitem[CS]{CS}
T. Coulhon \& L. Saloff-Coste, \emph{Isop\'erim\'etrie pour les
groupes et les vari\'et\'es}, Rev. Mat. Iberoamericana,
\textbf{9}~(1993), 293-314.


\bibitem[DGN1]{DGN1}
  D. Danielli, N. Garofalo \& D. M. Nhieu, \emph{Trace inequalities for Carnot-Carath\'eodory spaces and applications}, Ann. Sc. Norm. Sup. Pisa, Cl. Sci. (4), 2, \textbf{27}~(1998), 195-252.

\bibitem[DGN2]{DGN2}
\bysame, \emph{Non-doubling Ahlfors measures, Perimeter measures,
and the characterization of the trace spaces of Sobolev functions in
Carnot-Carath\'eodory spaces}, Memoirs of the Amer. Math. Soc.,
vol.182, 2006, p.119.







\bibitem[DGN3]{DGN4}
D. Danielli, N. Garofalo \& D. M. Nhieu, \emph{Notions of
convexity in Carnot groups}, Comm. Anal. and Geom., \textbf{11},
no.2,~(2003), 263-341.


\bibitem[DGN4]{DGN3}
\bysame, \emph{Sub-Riemannian calculus on hypersurfaces in Carnot
groups}, preprint, 2005.


\bibitem[DGN5]{DGN5}
\bysame, \emph{A notable family of entire intrinsic minimal graphs
in the Heisenberg group which are not perimeter minimizing},
preprint, 2005.





\bibitem[DGNP]{DGNP}
D. Danielli, N. Garofalo, D. M. Nhieu \& S. D. Pauls,
\emph{Instability of graphical strips and a positive answer to the
Bernstein problem in the Heisenberg group $\HH$}, preprint, 2006.

\bibitem[DG1]{DG1}
E. De Giorgi, \emph{Su una teoria generale della misura $(r-1)-$dimensionale in uno spazio a $r$ dimensioni}, Ann. Mat. Pura Appl., \textbf{36}~(1954), 191-213.

\bibitem[DG2]{DG2}
\bysame, \emph{Nuovi teoremi relativi alla misura $(r-1)$-dimensionale in uno spazio a $r$ dimensioni}, Ric. Mat., \textbf{4}~(1955), 95-113.

\bibitem[DG3]{DG3}
\bysame, \emph{Sulla propriet\`a isoperimetrica dell'ipersfera,
nella classe degli insiemi aventi frontiera orientata di misura
finita}, (Italian) Atti Accad. Naz. Lincei. Mem. Cl. Sci. Fis.
Mat. Nat. Sez. I (8), \textbf{5}~1958 33-44.


\bibitem[DCP]{DCP}
E. De Giorgi, F. Colombini \& L. C. Piccinini, \emph{Frontiere orientate di misura minima e questioni collegate}, Sc. Norm. Sup. Pisa, Cl. Scienze, Quaderni, 1972.

\bibitem[De1]{De1}
 M. Derridj , \emph{Un probl\'eme aux limites pour une classe d'op\'erateurs du second ordre hypoelliptiques}, Ann. Inst. Fourier, Grenoble, \textbf{21}, 4 (1971), 99-148.

\bibitem[De2]{De2}
\bysame , \emph{Sur un th\'eor\`eme de traces},
 Ann. Inst. Fourier, Grenoble, \textbf{22}, 2 (1972), 73-83.

\bibitem[E1]{E1}
P. Eberlein, \emph{Geometry of $2$-step nilpotent groups with a left
invariant metric}, Ann. Sci. \'Ecole Norm. Sup. (4) \textbf{27}~
(1994), no. 5, 611-660.

\bibitem[E2]{E2}
\bysame, \emph{Geometry of $2$-step nilpotent groups with a left invariant metric. II}, Trans. Amer. Math. Soc., \textbf{343}~ (1994), no. 2, 805-828.

\bibitem[E3]{E3}
\bysame, \emph{Geometry of nonpositively curved manifolds}, Chicago Lectures in Mathematics. University of Chicago Press, Chicago, IL, 1996.


\bibitem[Fe]{Fe}
H. Federer, \emph{Geometric Measure Theory}, Springer, 1969.

\bibitem[FMP]{FMP}
C. B. Figueroa, F. Mercuri \& R. H. L. Pedrosa, \emph{Invariant
surfaces of the Heisenberg groups}, Ann. Mat. Pura Appl. (4),
\textbf{177}~(1999), 173-194.

\bibitem[FR]{FR}
W. H. Fleming \& R. Rishel, \emph{An integral formula for total
gradient variation}, Arch. Math., \textbf{11}~(1960), 218-222.

\bibitem[F1]{F1}
G. B. Folland, \emph{A fundamental solution for a subelliptic
operator}, Bull. Amer. Math. Soc., \textbf{79}~(1973), 373-376.


\bibitem[F2]{F2}
\bysame, \emph{Subelliptic estimates and function spaces on
nilpotent Lie groups}, Ark. Math., \textbf{13}~(1975), 161-207.

\bibitem[F3]{F3}
\bysame, \emph{Harmonic Analysis in Phase Space}, Ann. Math.
Studies, Princeton Univ. Press, 1989.








\bibitem[FGW]{FGW}
  B. Franchi, S. Gallot \& R. Wheeden, \emph{Sobolev and isoperimetric inequalities for degenerate metrics}, Math. Ann., \textbf{300}~(1994), 557-571.



\bibitem[FSS1]{FSS1}
  B. Franchi, R. Serapioni \& F. Serra Cassano,
\emph{Meyers-Serrin type theorems and relaxation of variational integrals depending on vector fields}. Houston J. Math. \textbf{22}~(1996), no. 4, 859-890.

\bibitem[FSS2]{FSS2}
\bysame, \emph{Rectifiability and perimeter in the Heisenberg group}, Math. Ann., \textbf{321}~(2001) 3, 479-531.


\bibitem[FSS3]{FSS3}
\bysame, \emph{On the structure of finite perimeter sets in step
$2$ Carnot groups}, J. Geom. Anal., \textbf{13}~(2003), no. 3,
421-466.

\bibitem[FSS4]{FSS4}
\bysame, \emph{Regular hypersurfaces, intrinsic perimeter and
implicit function theorem in Carnot groups}, Comm. Anal. Geom.,
\textbf{11}~(2003), no. 5, 909-944.




\bibitem[GN]{GN}
  N. Garofalo \& D. M. Nhieu, \emph{Isoperimetric and Sobolev inequalities for Carnot-Carath\'eodory spaces and the existence of minimal surfaces}, Comm. Pure
Appl. Math., \textbf{49}~(1996), 1081-1144.


\bibitem[GP]{GP}
N. Garofalo \& S. D. Pauls, \emph{The Bernstein problem in the
Heisenberg group}, preprint, 2004.


\bibitem[GH]{GH}
M. Giaquinta \& S. Hildebrandt, \emph{Calculus of variations}, I
\& II. Grundlehren der Mathematischen Wissenschaften [Fundamental
Principles of Mathematical Sciences], 310, 311. Springer-Verlag,
Berlin, 1996.



\bibitem[Gro1]{Gro1}
M. Gromov, \emph{Carnot-Carath\'eodory spaces seen from within},
in Sub-Riemannian Geometry, Progress in Mathematics, vol. 144,
edited by Andr\'e Bella\"iche \& Jean-Jacques Risler,
Birkh\"auser, 1996.


\bibitem[Gro2]{Gro2}
\bysame, \emph{Metric Structures for Riemannian and Non-Riemannian Spaces}, Ed. by J. LaFontaine and P. Pansu, Birkh\"auser, 1998.

\bibitem[HP]{HP}
R. K. Hladky \& S. D. Pauls, \emph{Constant mean curvature
surfaces in sub-Riemannian geometry}, preprint, 2005.

\bibitem[H]{H}
  H. H\"ormander, \emph{Hypoelliptic second-order differential equations}, Acta Math., \textbf{119}~(1967), 147-171.





\bibitem[LM]{LM}
G. P. Leonardi \& S. Masnou, \emph{On the isoperimetric problem in
the Heisenberg group $\Hn$}, Ann. Mat. Pura Appl., (4)
\textbf{184}~(2005), no. 4, 533-553.


\bibitem[LR]{LR}
G. P. Leonardi \& S. Rigot, \emph{Isoperimetric sets on Carnot
groups}, Houston J. Math., \textbf{29}~(2003), no. 3, 609-637.


\bibitem[Ma]{Ma}
V. Magnani, \emph{Characteristic points, rectifiability and
perimeter measure on stratified groups}, preprint, 2004.


\bibitem[MaSC]{MaSC}
P. Maheux \& L. Saloff-Coste, \emph{Analyse sur les boules d'un op\'erateur sous-elliptique}, Math. Ann., \textbf{303}~(1995), 713-740.

\bibitem[MM]{MM}
U. Massari \& M. Miranda, \emph{Minimal Surfaces of Codimension
One}, Math. Studies \textbf{91}, North-Holland, 1984.

\bibitem[Mo1]{Mo1}
 R. Monti, \emph{Some properties of Carnot-Carathéodory balls in the Heisenberg
 group}, Atti Accad. Naz. Lincei Cl. Sci. Fis. Mat. Natur. Rend. Lincei
 (9), Mat. Appl.textbf{11}~(2000), no. 3, 155-167 (2001).

 \bibitem[Mo2]{Mo2}
 \bysame, \emph{Brunn-Minkowski and isoperimetric inequality in the Heisenberg group}, Ann. Acad. Sci. Fenn. Math., \textbf{28}~(2003), no. 1,
 99-109.

\bibitem[MoM]{MoM}
R. Monti \& D. Morbidelli, \emph{Isoperimetric inequality in the
Grushin plane}, J. Geom. Anal., \textbf{14}~(2004), no. 2,
355-368.

\bibitem[P1]{P}
P. Pansu, \emph{Une in\'egalit\'e isop\'erim\'etrique sur le
groupe de Heisenberg}, C. R. Acad. Sci. Paris S\'er. I Math.,
\textbf{295}~ (1982), no. 2, 127-130.


\bibitem[P2]{P2}
P. Pansu, \emph{An isoperimetric inequality on the Heisenberg
group}, Conference on differential geometry on homogeneous spaces
(Torino, 1983). Rend. Sem. Mat. Univ. Politec. Torino 1983,
Special Issue, 159-174 (1984).




\bibitem[Pa]{Pa}
S. Pauls, \emph{Minimal surfaces in the Heisenberg group}, Geom.
Dedicata, \textbf{104}~(2004), 201-231.


\bibitem[RR1]{RR}
M. Ritor\`e \& C. Rosales, \emph{Rotationally invariant
hypersurfaces with constant mean curvature in the Heisenberg group
$\Hn$}, preprint, April 2005.


\bibitem[RR2]{RR2}
\bysame, \emph{Area stationary surfaces in the Heisenberg group
$\HH$}, preprint, December 2005.


\bibitem[STh]{STh}
I. M. Singer \& J. A. Thorpe
\emph{Lecture notes on elementary topology and geometry},
Scott-Foresman, Glenview, Illinois, 1967.

\bibitem[Tr]{Tr}
  J. L. Troutman, \emph{Variational Calculus with Elementary Convexity},
  Springer-Verlag, 1983.



\bibitem[St]{St}
  E. M. Stein, \emph{Harmonic Analysis: Real Variable Methods, Orthogonality and Oscillatory Integrals}, Princeton Univ. Press, (1993).


\bibitem[To]{To}
P. Tomter, \emph{Constant mean curvature surfaces in the Heisenberg group}, \emph{Differential Geometry: Partial Differential Equations of Manifolds (Los Angeles, CA, 1990)}, p.485-495, Amer. Math. Soc., Providence, RI, 1993.



\bibitem[Va1]{Va1}
N. Th. Varopoulos, \emph{Analysis on nilpotent groups}, J. Funct.
Anal., \textbf{66}~(1986), no.3, 406-431.

\bibitem[Va2]{Va2}
\bysame, \emph{Analysis on Lie groups}, J. Funct. Anal., \textbf{67}~(1988), no.2, 346-410.

\bibitem[VSC]{VSC}
N.  Th. Varopoulos, L. Saloff-Coste \& T. Coulhon, \emph{Analysis and Geometry on Groups}, Cambridge U. Press, 1992.


\bibitem[Z]{Z}
 W. P. Ziemer, \emph{Weakly Differentiable Functions},
Springer-Verlag (1989).




\end{thebibliography}
\end{document}